\voffset=-0.75in
\magnification=1200
\vsize=8.00in\hsize = 7.00 true in\hoffset= -.25in
\tolerance = 10000
$\ $ \vskip 1.00in
\centerline{``Patterns of Compact Cardinals''}
\vskip .5in
\centerline{by}
\vskip .5in
\centerline{Arthur W. Apter${}^*$}
\centerline{Department of Mathematics}
\centerline{Baruch College of CUNY}
\centerline{New York, New York 10010}
\vskip .25in
\centerline{November 14, 1996}
\centerline{(Revised February 20, 1997)}

\vskip .50in
\noindent Abstract: We show relative to strong hypotheses
that patterns of compact cardinals in the universe,
where a compact cardinal is one which is either
strongly compact or supercompact, can be
virtually arbitrary.  Specifically, we prove if
$V \models ``$ZFC + $\Omega$ is the least inaccessible
limit of measurable limits of supercompact cardinals +
$f : \Omega \to 2$ is a function'', then there is
a partial ordering $P \in V$ so that for
$\overline V = V^P$, $\overline V_\Omega \models
``$ZFC + There is a proper class of compact cardinals +
If $f(\alpha) = 0$, then the
$\alpha^{{\hbox{\rm th}}}$
compact cardinal isn't supercompact +
If $f(\alpha) = 1$, then the
$\alpha^{{\hbox{\rm th}}}$
compact cardinal is supercompact''.  We then prove a
generalized version of this theorem assuming $\kappa$ is
a supercompact limit of supercompact cardinals and
$f : \kappa \to 2$ is a function, and we derive as
corollaries of the generalized version of the theorem
the consistency of the least measurable limit of
supercompact cardinals being the same as the least
measurable limit of non-supercompact strongly compact
cardinals and the consistency of the least supercompact
cardinal being a limit of strongly compact cardinals. \hfil\break


cardinal. \hfil\break

\noindent ${}^*$The author wishes to thank Joel Hamkins
for helpful conversations on the subject matter
of this paper.  In additon, the author wishes to thank
the referee, both for promptly reviewing the paper and for
making helpful comments and suggestions which improved
the presentation of the material contained herein.
\eject
\catcode`\@=11
\font\tensmc=cmcsc10      
\def\smc{\tensmc}

\def\hcorrection#1{\advance\hoffset by #1 }
\def\vcorrection#1{\advance\voffset by #1 }
\def\wlog#1{}
\newif\iftitle@
\outer\def\title{\title@true\vglue 24\p@ plus 12\p@ minus 12\p@
   \bgroup\let\\=\cr\tabskip\centering
   \halign to \hsize\bgroup\tenbf\hfill\ignorespaces##\unskip\hfill\cr}
\def\endtitle{\cr\egroup\egroup\vglue 18\p@ plus 12\p@ minus 6\p@}
\outer\def\author{\iftitle@\vglue -18\p@ plus -12\p@ minus -6\p@\fi\vglue
    12\p@ plus 6\p@ minus 3\p@\bgroup\let\\=\cr\tabskip\centering
    \halign to \hsize\bgroup\smc\hfill\ignorespaces##\unskip\hfill\cr}
\def\endauthor{\cr\egroup\egroup\vglue 18\p@ plus 12\p@ minus 6\p@}
\outer\def\heading{\bigbreak\bgroup\let\\=\cr\tabskip\centering
    \halign to \hsize\bgroup\smc\hfill\ignorespaces##\unskip\hfill\cr}
\def\endheading{\cr\egroup\egroup\nobreak\medskip}

\outer\def\proclaim#1{\medbreak\noindent\smc\ignorespaces
    #1\unskip.\enspace\sl\ignorespaces}
\outer\def\endproclaim{\par\ifdim\lastskip<\medskipamount\removelastskip
  \penalty 55 \fi\medskip\rm}
\outer\def\demo#1{\par\ifdim\lastskip<\smallskipamount\removelastskip
    \smallskip\fi\noindent{\smc\ignorespaces#1\unskip:\enspace}\rm
      \ignorespaces}

\newcount\footmarkcount@
\footmarkcount@=1
\def\makefootnote@#1#2{\insert\footins{\interlinepenalty=100
  \splittopskip=\ht\strutbox \splitmaxdepth=\dp\strutbox
  \floatingpenalty=\@MM
  \leftskip=\z@\rightskip=\z@\spaceskip=\z@\xspaceskip=\z@
  \noindent{#1}\footstrut\rm\ignorespaces #2\strut}}
\def\footnote{\let\@sf=\empty\ifhmode\edef\@sf{\spacefactor
   =\the\spacefactor}\/\fi\futurelet\next\footnote@}
\def\footnote@{\ifx"\next\let\next\footnote@@\else
    \let\next\footnote@@@\fi\next}
\def\footnote@@"#1"#2{#1\@sf\relax\makefootnote@{#1}{#2}}
\def\footnote@@@#1{$^{\number\footmarkcount@}$\makefootnote@
   {$^{\number\footmarkcount@}$}{#1}\global\advance\footmarkcount@ by 1 }

\hyphenation{man-u-script man-u-scripts ap-pen-dix ap-pen-di-ces}
\hyphenation{data-base data-bases}
\ifx\amstexloaded@\relax\catcode`\@=13
  \endinput\else\let\amstexloaded@=\relax\fi
\newlinechar=`\^^J
\def\eat@#1{}
\def\Space@.{\futurelet\Space@\relax}
\Space@. %
\newhelp\athelp@
{Only certain combinations beginning with @ make sense to me.^^J
Perhaps you wanted \string\@\space for a printed @?^^J
I've ignored the character or group after @.}
\def\futureletnextat@{\futurelet\next\at@}
{\catcode`\@=\active
\lccode`\Z=`\@ \lowercase
{\gdef@{\expandafter\csname futureletnextatZ\endcsname}
\expandafter\gdef\csname atZ\endcsname
   {\ifcat\noexpand\next a\def\next{\csname atZZ\endcsname}\else
   \ifcat\noexpand\next0\def\next{\csname atZZ\endcsname}\else
    \def\next{\csname atZZZ\endcsname}\fi\fi\next}
\expandafter\gdef\csname atZZ\endcsname#1{\expandafter
   \ifx\csname #1Zat\endcsname\relax\def\next
     {\errhelp\expandafter=\csname athelpZ\endcsname
      \errmessage{Invalid use of \string@}}\else
       \def\next{\csname #1Zat\endcsname}\fi\next}
\expandafter\gdef\csname atZZZ\endcsname#1{\errhelp
    \expandafter=\csname athelpZ\endcsname
      \errmessage{Invalid use of \string@}}}}
\def\atdef@#1{\expandafter\def\csname #1@at\endcsname}
\newhelp\defahelp@{If you typed \string\define\space cs instead of
\string\define\string\cs\space^^J
I've substituted an inaccessible control sequence so that your^^J
definition will be completed without mixing me up too badly.^^J
If you typed \string\define{\string\cs} the inaccessible control sequence^^J
was defined to be \string\cs, and the rest of your^^J
definition appears as input.}
\newhelp\defbhelp@{I've ignored your definition, because it might^^J
conflict with other uses that are important to me.}
\def\define{\futurelet\next\define@}
\def\define@{\ifcat\noexpand\next\relax
  \def\next{\define@@}%
  \else\errhelp=\defahelp@
  \errmessage{\string\define\space must be followed by a control
     sequence}\def\next{\def\garbage@}\fi\next}
\def\undefined@{}
\def\preloaded@{}
\def\define@@#1{\ifx#1\relax\errhelp=\defbhelp@
   \errmessage{\string#1\space is already defined}\def\next{\def\garbage@}%
   \else\expandafter\ifx\csname\expandafter\eat@\string
	 #1@\endcsname\undefined@\errhelp=\defbhelp@
   \errmessage{\string#1\space can't be defined}\def\next{\def\garbage@}%
   \else\expandafter\ifx\csname\expandafter\eat@\string#1\endcsname\relax
     \def\next{\def#1}\else\errhelp=\defbhelp@
     \errmessage{\string#1\space is already defined}\def\next{\def\garbage@}%
      \fi\fi\fi\next}
\def\famzero{\fam\z@}

\def\lim{\mathop{\famzero lim}}

\def\sup{\mathop{\famzero sup}}

\def\textfont@#1#2{\def#1{\relax\ifmmode
    \errmessage{Use \string#1\space only in text}\else#2\fi}}
\textfont@\rm\tenrm
\textfont@\it\tenit
\textfont@\sl\tensl
\textfont@\bf\tenbf
\textfont@\smc\tensmc
\let\ic@=\/
\def\/{\unskip\ic@}
\def\textfonti{\the\textfont1 }
\def\t#1#2{{\edef\next{\the\font}\textfonti\accent"7F \next#1#2}}
\let\B=\=
\let\D=\.
\def~{\unskip\nobreak\ \ignorespaces}
{\catcode`\@=\active
\gdef\@{\char'100 }}
\atdef@-{\leavevmode\futurelet\next\athyph@}
\def\athyph@{\ifx\next-\let\next=\athyph@@
  \else\let\next=\athyph@@@\fi\next}
\def\athyph@@@{\hbox{-}}
\def\athyph@@#1{\futurelet\next\athyph@@@@}
\def\athyph@@@@{\if\next-\def\next##1{\hbox{---}}\else
    \def\next{\hbox{--}}\fi\next}
\def\.{.\spacefactor=\@m}
\atdef@.{\null.}
\atdef@,{\null,}
\atdef@;{\null;}
\atdef@:{\null:}
\atdef@?{\null?}
\atdef@!{\null!}
\def\srdr@{\thinspace}
\def\drsr@{\kern.02778em}
\def\sldl@{\kern.02778em}
\def\dlsl@{\thinspace}
\atdef@"{\unskip\futurelet\next\atqq@}
\def\atqq@{\ifx\next\Space@\def\next. {\atqq@@}\else
	 \def\next.{\atqq@@}\fi\next.}
\def\atqq@@{\futurelet\next\atqq@@@}
\def\atqq@@@{\ifx\next`\def\next`{\atqql@}\else\def\next'{\atqqr@}\fi\next}
\def\atqql@{\futurelet\next\atqql@@}
\def\atqql@@{\ifx\next`\def\next`{\sldl@``}\else\def\next{\dlsl@`}\fi\next}
\def\atqqr@{\futurelet\next\atqqr@@}
\def\atqqr@@{\ifx\next'\def\next'{\srdr@''}\else\def\next{\drsr@'}\fi\next}
\def\flushpar{\par\noindent}
\def\textfontii{\the\textfont2 }
\def\{{\relax\ifmmode\lbrace\else
    {\textfontii f}\spacefactor=\@m\fi}
\def\}{\relax\ifmmode\rbrace\else
    \let\@sf=\empty\ifhmode\edef\@sf{\spacefactor=\the\spacefactor}\fi
      {\textfontii g}\@sf\relax\fi}
\def\nonhmodeerr@#1{\errmessage
     {\string#1\space allowed only within text}}
\def\linebreak{\relax\ifhmode\unskip\break\else
    \nonhmodeerr@\linebreak\fi}
\def\allowlinebreak{\relax
   \ifhmode\allowbreak\else\nonhmodeerr@\allowlinebreak\fi}
\newskip\saveskip@
\def\nolinebreak{\relax\ifhmode\saveskip@=\lastskip\unskip
  \nobreak\ifdim\saveskip@>\z@\hskip\saveskip@\fi
   \else\nonhmodeerr@\nolinebreak\fi}
\def\newline{\relax\ifhmode\null\hfil\break
    \else\nonhmodeerr@\newline\fi}
\def\nonmathaerr@#1{\errmessage
     {\string#1\space is not allowed in display math mode}}
\def\nonmathberr@#1{\errmessage{\string#1\space is allowed only in math mode}}
\def\mathbreak{\relax\ifmmode\ifinner\break\else
   \nonmathaerr@\mathbreak\fi\else\nonmathberr@\mathbreak\fi}
\def\nomathbreak{\relax\ifmmode\ifinner\nobreak\else
    \nonmathaerr@\nomathbreak\fi\else\nonmathberr@\nomathbreak\fi}
\def\allowmathbreak{\relax\ifmmode\ifinner\allowbreak\else
     \nonmathaerr@\allowmathbreak\fi\else\nonmathberr@\allowmathbreak\fi}
\def\pagebreak{\relax\ifmmode
   \ifinner\errmessage{\string\pagebreak\space
     not allowed in non-display math mode}\else\postdisplaypenalty-\@M\fi
   \else\ifvmode\penalty-\@M\else\edef\spacefactor@
       {\spacefactor=\the\spacefactor}\vadjust{\penalty-\@M}\spacefactor@
	\relax\fi\fi}
\def\nopagebreak{\relax\ifmmode
     \ifinner\errmessage{\string\nopagebreak\space
    not allowed in non-display math mode}\else\postdisplaypenalty\@M\fi
    \else\ifvmode\nobreak\else\edef\spacefactor@
	{\spacefactor=\the\spacefactor}\vadjust{\penalty\@M}\spacefactor@
	 \relax\fi\fi}
\def\newpage{\relax\ifvmode\vfill\penalty-\@M\else\nonvmodeerr@\newpage\fi}
\def\nonvmodeerr@#1{\errmessage
    {\string#1\space is allowed only between paragraphs}}
\def\smallpagebreak{\relax\ifvmode\smallbreak
      \else\nonvmodeerr@\smallpagebreak\fi}
\def\medpagebreak{\relax\ifvmode\medbreak
       \else\nonvmodeerr@\medpagebreak\fi}
\def\bigpagebreak{\relax\ifvmode\bigbreak
      \else\nonvmodeerr@\bigpagebreak\fi}
\newdimen\captionwidth@
\captionwidth@=\hsize
\advance\captionwidth@ by -1.5in
\def\caption#1{}
\def\topspace#1{\gdef\thespace@{#1}\ifvmode\def\next
    {\futurelet\next\topspace@}\else\def\next{\nonvmodeerr@\topspace}\fi\next}
\def\topspace@{\ifx\next\Space@\def\next. {\futurelet\next\topspace@@}\else
     \def\next.{\futurelet\next\topspace@@}\fi\next.}
\def\topspace@@{\ifx\next\caption\let\next\topspace@@@\else
    \let\next\topspace@@@@\fi\next}
 \def\topspace@@@@{\topinsert\vbox to
       \thespace@{}\endinsert}
\def\topspace@@@\caption#1{\topinsert\vbox to
    \thespace@{}\nobreak
      \smallskip
    \setbox\z@=\hbox{\noindent\ignorespaces#1\unskip}%
   \ifdim\wd\z@>\captionwidth@
   \centerline{\vbox{\hsize=\captionwidth@\noindent\ignorespaces#1\unskip}}%
   \else\centerline{\box\z@}\fi\endinsert}
\def\midspace#1{\gdef\thespace@{#1}\ifvmode\def\next
    {\futurelet\next\midspace@}\else\def\next{\nonvmodeerr@\midspace}\fi\next}
\def\midspace@{\ifx\next\Space@\def\next. {\futurelet\next\midspace@@}\else
     \def\next.{\futurelet\next\midspace@@}\fi\next.}
\def\midspace@@{\ifx\next\caption\let\next\midspace@@@\else
    \let\next\midspace@@@@\fi\next}
 \def\midspace@@@@{\midinsert\vbox to
       \thespace@{}\endinsert}
\def\midspace@@@\caption#1{\midinsert\vbox to
    \thespace@{}\nobreak
      \smallskip
      \setbox\z@=\hbox{\noindent\ignorespaces#1\unskip}%
      \ifdim\wd\z@>\captionwidth@
    \centerline{\vbox{\hsize=\captionwidth@\noindent\ignorespaces#1\unskip}}%
    \else\centerline{\box\z@}\fi\endinsert}
\mathchardef\prime@="0230
\def\prime{{{}\prime@{}}}
\def\prim@s{\prime@\futurelet\next\pr@m@s}

\def\,{\relax\ifmmode\mskip\thinmuskip\else\thinspace\fi}
\def\!{\relax\ifmmode\mskip-\thinmuskip\else\negthinspace\fi}
\def\frac#1#2{{#1\over#2}}

\def\:{\nobreak\hskip.1111em{:}\hskip.3333em plus .0555em\relax}
\def\intic@{\mathchoice{\hskip5\p@}{\hskip4\p@}{\hskip4\p@}{\hskip4\p@}}
\def\negintic@
 {\mathchoice{\hskip-5\p@}{\hskip-4\p@}{\hskip-4\p@}{\hskip-4\p@}}
\def\intkern@{\mathchoice{\!\!\!}{\!\!}{\!\!}{\!\!}}
\def\intdots@{\mathchoice{\cdots}{{\cdotp}\mkern1.5mu
    {\cdotp}\mkern1.5mu{\cdotp}}{{\cdotp}\mkern1mu{\cdotp}\mkern1mu
      {\cdotp}}{{\cdotp}\mkern1mu{\cdotp}\mkern1mu{\cdotp}}}
\newcount\intno@
\def\iint{\intno@=\tw@\futurelet\next\ints@}
\def\iiint{\intno@=\thr@@\futurelet\next\ints@}
\def\iiiint{\intno@=4 \futurelet\next\ints@}
\def\idotsint{\intno@=\z@\futurelet\next\ints@}
\def\ints@{\findlimits@\ints@@}
\newif\iflimtoken@
\newif\iflimits@
\def\findlimits@{\limtoken@false\limits@false\ifx\next\limits
 \limtoken@true\limits@true\else\ifx\next\nolimits\limtoken@true\limits@false
    \fi\fi}
\def\multintlimits@{\intop\ifnum\intno@=\z@\intdots@
  \else\intkern@\fi
    \ifnum\intno@>\tw@\intop\intkern@\fi
     \ifnum\intno@>\thr@@\intop\intkern@\fi\intop}
\def\multint@{\int\ifnum\intno@=\z@\intdots@\else\intkern@\fi
   \ifnum\intno@>\tw@\int\intkern@\fi
    \ifnum\intno@>\thr@@\int\intkern@\fi\int}
\def\ints@@{\iflimtoken@\def\ints@@@{\iflimits@
   \negintic@\mathop{\intic@\multintlimits@}\limits\else
    \multint@\nolimits\fi\eat@}\else
     \def\ints@@@{\multint@\nolimits}\fi\ints@@@}
\def\Sb{_\bgroup\vspace@
	\baselineskip=\fontdimen10 \scriptfont\tw@
	\advance\baselineskip by \fontdimen12 \scriptfont\tw@
	\lineskip=\thr@@\fontdimen8 \scriptfont\thr@@
	\lineskiplimit=\thr@@\fontdimen8 \scriptfont\thr@@
	\Let@\vbox\bgroup\halign\bgroup \hfil$\scriptstyle
	    {##}$\hfil\cr}
\def\endSb{\crcr\egroup\egroup\egroup}
\def\Sp{^\bgroup\vspace@
	\baselineskip=\fontdimen10 \scriptfont\tw@
	\advance\baselineskip by \fontdimen12 \scriptfont\tw@
	\lineskip=\thr@@\fontdimen8 \scriptfont\thr@@
	\lineskiplimit=\thr@@\fontdimen8 \scriptfont\thr@@
	\Let@\vbox\bgroup\halign\bgroup \hfil$\scriptstyle
	    {##}$\hfil\cr}
\def\endSp{\crcr\egroup\egroup\egroup}
\def\Let@{\relax\iffalse{\fi\let\\=\cr\iffalse}\fi}
\def\vspace@{\def\vspace##1{\noalign{\vskip##1 }}}
\def\aligned{\,\vcenter\bgroup\vspace@\Let@\openup\jot\m@th\ialign
  \bgroup \strut\hfil$\displaystyle{##}$&$\displaystyle{{}##}$\hfil\crcr}
\def\endaligned{\crcr\egroup\egroup}
\def\matrix{\,\vcenter\bgroup\Let@\vspace@
    \normalbaselines
  \m@th\ialign\bgroup\hfil$##$\hfil&&\quad\hfil$##$\hfil\crcr
    \mathstrut\crcr\noalign{\kern-\baselineskip}}
\def\endmatrix{\crcr\mathstrut\crcr\noalign{\kern-\baselineskip}\egroup
		\egroup\,}
\newtoks\hashtoks@
\hashtoks@={#}
\def\format{\crcr\egroup\iffalse{\fi\ifnum`}=0 \fi\format@}
\def\format@#1\\{\def\preamble@{#1}%
  \def\c{\hfil$\the\hashtoks@$\hfil}%
  \def\r{\hfil$\the\hashtoks@$}%
  \def\l{$\the\hashtoks@$\hfil}%
  \setbox\z@=\hbox{\xdef\Preamble@{\preamble@}}\ifnum`{=0 \fi\iffalse}\fi
   \ialign\bgroup\span\Preamble@\crcr}

\def\cases{\left\{\,\vcenter\bgroup\vspace@
     \normalbaselines\openup\jot\m@th
       \Let@\ialign\bgroup$##$\hfil&\quad$##$\hfil\crcr
      \mathstrut\crcr\noalign{\kern-\baselineskip}}

\newif\iftagsleft@
\tagsleft@true
\def\TagsOnRight{\global\tagsleft@false}
\def\tag#1$${\iftagsleft@\leqno\else\eqno\fi
 \hbox{\def\pagebreak{\global\postdisplaypenalty-\@M}%
 \def\nopagebreak{\global\postdisplaypenalty\@M}\rm(#1\unskip)}%
  $$\postdisplaypenalty\z@\ignorespaces}
\interdisplaylinepenalty=\@M
\def\allowdisplaybreak@{\def\allowdisplaybreak{\noalign{\allowbreak}}}
\def\displaybreak@{\def\displaybreak{\noalign{\break}}}
\def\align#1\endalign{\def\tag{&}\vspace@\allowdisplaybreak@\displaybreak@
  \iftagsleft@\lalign@#1\endalign\else
   \ralign@#1\endalign\fi}
\def\ralign@#1\endalign{\displ@y\Let@\tabskip\centering\halign to\displaywidth
     {\hfil$\displaystyle{##}$\tabskip=\z@&$\displaystyle{{}##}$\hfil
       \tabskip=\centering&\llap{\hbox{(\rm##\unskip)}}\tabskip\z@\crcr
	     #1\crcr}}
\def\lalign@
 #1\endalign{\displ@y\Let@\tabskip\centering\halign to \displaywidth
   {\hfil$\displaystyle{##}$\tabskip=\z@&$\displaystyle{{}##}$\hfil
   \tabskip=\centering&\kern-\displaywidth
	\rlap{\hbox{(\rm##\unskip)}}\tabskip=\displaywidth\crcr
	       #1\crcr}}
\def\overrightarrow{\mathpalette\overrightarrow@}
\def\overrightarrow@#1#2{\vbox{\ialign{$##$\cr
    #1{-}\mkern-6mu\cleaders\hbox{$#1\mkern-2mu{-}\mkern-2mu$}\hfill
     \mkern-6mu{\to}\cr
     \noalign{\kern -1\p@\nointerlineskip}
     \hfil#1#2\hfil\cr}}}
\def\overleftarrow{\mathpalette\overleftarrow@}
\def\overleftarrow@#1#2{\vbox{\ialign{$##$\cr
     #1{\leftarrow}\mkern-6mu\cleaders\hbox{$#1\mkern-2mu{-}\mkern-2mu$}\hfill
      \mkern-6mu{-}\cr
     \noalign{\kern -1\p@\nointerlineskip}
     \hfil#1#2\hfil\cr}}}
\def\overleftrightarrow{\mathpalette\overleftrightarrow@}
\def\overleftrightarrow@#1#2{\vbox{\ialign{$##$\cr
     #1{\leftarrow}\mkern-6mu\cleaders\hbox{$#1\mkern-2mu{-}\mkern-2mu$}\hfill
       \mkern-6mu{\to}\cr
    \noalign{\kern -1\p@\nointerlineskip}
      \hfil#1#2\hfil\cr}}}
\def\underrightarrow{\mathpalette\underrightarrow@}
\def\underrightarrow@#1#2{\vtop{\ialign{$##$\cr
    \hfil#1#2\hfil\cr
     \noalign{\kern -1\p@\nointerlineskip}
    #1{-}\mkern-6mu\cleaders\hbox{$#1\mkern-2mu{-}\mkern-2mu$}\hfill
     \mkern-6mu{\to}\cr}}}
\def\underleftarrow{\mathpalette\underleftarrow@}
\def\underleftarrow@#1#2{\vtop{\ialign{$##$\cr
     \hfil#1#2\hfil\cr
     \noalign{\kern -1\p@\nointerlineskip}
     #1{\leftarrow}\mkern-6mu\cleaders\hbox{$#1\mkern-2mu{-}\mkern-2mu$}\hfill
      \mkern-6mu{-}\cr}}}
\def\underleftrightarrow{\mathpalette\underleftrightarrow@}
\def\underleftrightarrow@#1#2{\vtop{\ialign{$##$\cr
      \hfil#1#2\hfil\cr
    \noalign{\kern -1\p@\nointerlineskip}
     #1{\leftarrow}\mkern-6mu\cleaders\hbox{$#1\mkern-2mu{-}\mkern-2mu$}\hfill
       \mkern-6mu{\to}\cr}}}
\def\sqrt#1{\radical"270370 {#1}}
\def\dots{\relax\ifmmode\let\next=\ldots\else\let\next=\tdots@\fi\next}
\def\tdots@{\unskip\ \tdots@@}
\def\tdots@@{\futurelet\next\tdots@@@}
\def\tdots@@@{$\mathinner{\ldotp\ldotp\ldotp}\,
   \ifx\next,$\else
   \ifx\next.\,$\else
   \ifx\next;\,$\else
   \ifx\next:\,$\else
   \ifx\next?\,$\else
   \ifx\next!\,$\else
   $ \fi\fi\fi\fi\fi\fi}
\def\text{\relax\ifmmode\let\next=\text@\else\let\next=\text@@\fi\next}
\def\text@@#1{\hbox{#1}}
\def\text@#1{\mathchoice
 {\hbox{\everymath{\displaystyle}\def\textfonti{\the\textfont1 }%
    \def\textfontii{\the\textfont2 }\textdef@@ T#1}}
 {\hbox{\everymath{\textstyle}\def\textfonti{\the\textfont1 }%
    \def\textfontii{\the\textfont2 }\textdef@@ T#1}}
 {\hbox{\everymath{\scriptstyle}\def\textfonti{\the\scriptfont1 }%
   \def\textfontii{\the\scriptfont2 }\textdef@@ S\rm#1}}
 {\hbox{\everymath{\scriptscriptstyle}\def\textfonti{\the\scriptscriptfont1 }%
   \def\textfontii{\the\scriptscriptfont2 }\textdef@@ s\rm#1}}}
\def\textdef@@#1{\textdef@#1\rm \textdef@#1\bf
   \textdef@#1\sl \textdef@#1\it}

\def\textdef@#1#2{\def\next{\csname\expandafter\eat@\string#2fam\endcsname}%
\if S#1\edef#2{\the\scriptfont\next\relax}%
 \else\if s#1\edef#2{\the\scriptscriptfont\next\relax}%
 \else\edef#2{\the\textfont\next\relax}\fi\fi}
\scriptfont\itfam=\tenit \scriptscriptfont\itfam=\tenit
\scriptfont\slfam=\tensl \scriptscriptfont\slfam=\tensl
\mathcode`\0="0030
\mathcode`\1="0031
\mathcode`\2="0032
\mathcode`\3="0033
\mathcode`\4="0034
\mathcode`\5="0035
\mathcode`\6="0036
\mathcode`\7="0037
\mathcode`\8="0038
\mathcode`\9="0039
\def\Cal{\relax\ifmmode\let\next=\Cal@\else
     \def\next{\errmessage{Use \string\Cal\space only in math mode}}\fi\next}
\def\Cal@#1{{\fam2 #1}}
\def\bold{\relax\ifmmode\let\next=\bold@\else
   \def\next{\errmessage{Use \string\bold\space only in math
      mode}}\fi\next}\def\bold@#1{{\fam\bffam #1}}
\mathchardef\Gamma="0000
\mathchardef\Delta="0001
\mathchardef\Theta="0002
\mathchardef\Lambda="0003
\mathchardef\Xi="0004
\mathchardef\Pi="0005
\mathchardef\Sigma="0006
\mathchardef\Upsilon="0007
\mathchardef\Phi="0008
\mathchardef\Psi="0009
\mathchardef\Omega="000A
\mathchardef\varGamma="0100
\mathchardef\varDelta="0101
\mathchardef\varTheta="0102
\mathchardef\varLambda="0103
\mathchardef\varXi="0104
\mathchardef\varPi="0105
\mathchardef\varSigma="0106
\mathchardef\varUpsilon="0107
\mathchardef\varPhi="0108
\mathchardef\varPsi="0109
\mathchardef\varOmega="010A
\font\dummyft@=dummy
\fontdimen1 \dummyft@=\z@
\fontdimen2 \dummyft@=\z@
\fontdimen3 \dummyft@=\z@
\fontdimen4 \dummyft@=\z@
\fontdimen5 \dummyft@=\z@
\fontdimen6 \dummyft@=\z@
\fontdimen7 \dummyft@=\z@
\fontdimen8 \dummyft@=\z@
\fontdimen9 \dummyft@=\z@
\fontdimen10 \dummyft@=\z@
\fontdimen11 \dummyft@=\z@
\fontdimen12 \dummyft@=\z@
\fontdimen13 \dummyft@=\z@
\fontdimen14 \dummyft@=\z@
\fontdimen15 \dummyft@=\z@
\fontdimen16 \dummyft@=\z@
\fontdimen17 \dummyft@=\z@
\fontdimen18 \dummyft@=\z@
\fontdimen19 \dummyft@=\z@
\fontdimen20 \dummyft@=\z@
\fontdimen21 \dummyft@=\z@
\fontdimen22 \dummyft@=\z@
\def\fontlist@{\\{\tenrm}\\{\sevenrm}\\{\fiverm}\\{\teni}\\{\seveni}%
 \\{\fivei}\\{\tensy}\\{\sevensy}\\{\fivesy}\\{\tenex}\\{\tenbf}\\{\sevenbf}%
 \\{\fivebf}\\{\tensl}\\{\tenit}\\{\tensmc}}
\def\dodummy@{{\def\\##1{\global\let##1=\dummyft@}\fontlist@}}
\newif\ifsyntax@
\newcount\countxviii@
\def\newtoks@{\alloc@5\toks\toksdef\@cclvi}
\def\nopages@{\output={\setbox\z@=\box\@cclv \deadcycles=\z@}\newtoks@\output}
\def\syntax{\syntax@true\dodummy@\countxviii@=\count18
\loop \ifnum\countxviii@ > \z@ \textfont\countxviii@=\dummyft@
   \scriptfont\countxviii@=\dummyft@ \scriptscriptfont\countxviii@=\dummyft@
     \advance\countxviii@ by-\@ne\repeat
\dummyft@\tracinglostchars=\z@
  \nopages@\frenchspacing\hbadness=\@M}
\def\magstep#1{\ifcase#1 1000\or
 1200\or 1440\or 1728\or 2074\or 2488\or
 \errmessage{\string\magstep\space only works up to 5}\fi\relax}
{\lccode`\2=`\p \lccode`\3=`\t
 \lowercase{\gdef\tru@#123{#1truept}}}

\def\scaletype#1{\mag=#1\relax
 \hsize=\expandafter\tru@\the\hsize
 \vsize=\expandafter\tru@\the\vsize
 \dimen\footins=\expandafter\tru@\the\dimen\footins}

\def\scalefont#1#2\andcallit#3{\edef\font@{\the\font}#1\font#3=
  \fontname\font\space scaled #2\relax\font@}
\def\Mag@#1#2{\ifdim#1<1pt\multiply#1 #2\relax\divide#1 1000 \else
  \ifdim#1<10pt\divide#1 10 \multiply#1 #2\relax\divide#1 100\else
  \divide#1 100 \multiply#1 #2\relax\divide#1 10 \fi\fi}
\def\scalelinespacing#1{\Mag@\baselineskip{#1}\Mag@\lineskip{#1}%
  \Mag@\lineskiplimit{#1}}
\def\wlog#1{\immediate\write-1{#1}}
\catcode`\@=\active

\catcode`@=11
\def\binrel@#1{\setbox\z@\hbox{\thinmuskip0mu
\medmuskip\m@ne mu\thickmuskip\@ne mu$#1\m@th$}%
\setbox\@ne\hbox{\thinmuskip0mu\medmuskip\m@ne mu\thickmuskip
\@ne mu${}#1{}\m@th$}%
\setbox\tw@\hbox{\hskip\wd\@ne\hskip-\wd\z@}}
\def\overset#1\to#2{\binrel@{#2}\ifdim\wd\tw@<\z@
\mathbin{\mathop{\kern\z@#2}\limits^{#1}}\else\ifdim\wd\tw@>\z@
\mathrel{\mathop{\kern\z@#2}\limits^{#1}}\else
{\mathop{\kern\z@#2}\limits^{#1}}{}\fi\fi}
\def\underset#1\to#2{\binrel@{#2}\ifdim\wd\tw@<\z@
\mathbin{\mathop{\kern\z@#2}\limits_{#1}}\else\ifdim\wd\tw@>\z@
\mathrel{\mathop{\kern\z@#2}\limits_{#1}}\else
{\mathop{\kern\z@#2}\limits_{#1}}{}\fi\fi}
\def\circle#1{\leavevmode\setbox0=\hbox{h}\dimen@=\ht0
\advance\dimen@ by-1ex\rlap{\raise1.5\dimen@\hbox{\char'27}}#1}
\def\sqr#1#2{{\vcenter{\hrule height.#2pt
     \hbox{\vrule width.#2pt height#1pt \kern#1pt
       \vrule width.#2pt}
     \hrule height.#2pt}}}
\def\square{\mathchoice\sqr34\sqr34\sqr{2.1}3\sqr{1.5}3}
\def\force{\hbox{$\|\hskip-2pt\hbox{--}$\hskip2pt}}
 
\catcode`@=\active
 
\catcode`\@=11
\def\bold{\relaxnext@\ifmmode\let\next\bold@\else
 \def\next{\Err@{Use \string\bold\space only in math mode}}\fi\next}
\def\bold@#1{{\bold@@{#1}}}
\def\bold@@#1{\fam\bffam#1}
\def\hexnumber@#1{\ifnum#1<10 \number#1\else
 \ifnum#1=10 A\else\ifnum#1=11 B\else\ifnum#1=12 C\else
 \ifnum#1=13 D\else\ifnum#1=14 E\else\ifnum#1=15 F\fi\fi\fi\fi\fi\fi\fi}
\def\bffam@{\hexnumber@\bffam}
 
 
\font\tenmsx=msam10
\font\sevenmsx=msam7
\font\fivemsx=msam5
\font\tenmsy=msbm10
\font\sevenmsy=msbm7
\font\fivemsy=msbm7
 
\newfam\msxfam
\newfam\msyfam
\textfont\msxfam=\tenmsx
\scriptfont\msxfam=\sevenmsx
\scriptscriptfont\msxfam=\fivemsx
\textfont\msyfam=\tenmsy
\scriptfont\msyfam=\sevenmsy
\scriptscriptfont\msyfam=\fivemsy
\def\msx@{\hexnumber@\msxfam}
\def\msy@{\hexnumber@\msyfam}
\mathchardef\boxdot="2\msx@00
\mathchardef\boxplus="2\msx@01
\mathchardef\boxtimes="2\msx@02
\mathchardef\square="0\msx@03
\mathchardef\blacksquare="0\msx@04
\mathchardef\centerdot="2\msx@05
\mathchardef\lozenge="0\msx@06
\mathchardef\blacklozenge="0\msx@07
\mathchardef\circlearrowright="3\msx@08
\mathchardef\circlearrowleft="3\msx@09
\mathchardef\rightleftharpoons="3\msx@0A
\mathchardef\leftrightharpoons="3\msx@0B
\mathchardef\boxminus="2\msx@0C
\mathchardef\Vdash="3\msx@0D
\mathchardef\Vvdash="3\msx@0E
\mathchardef\vDash="3\msx@0F
\mathchardef\twoheadrightarrow="3\msx@10
\mathchardef\twoheadleftarrow="3\msx@11
\mathchardef\leftleftarrows="3\msx@12
\mathchardef\rightrightarrows="3\msx@13
\mathchardef\upuparrows="3\msx@14
\mathchardef\downdownarrows="3\msx@15
\mathchardef\upharpoonright="3\msx@16

\mathchardef\downharpoonright="3\msx@17
\mathchardef\upharpoonleft="3\msx@18
\mathchardef\downharpoonleft="3\msx@19
\mathchardef\rightarrowtail="3\msx@1A
\mathchardef\leftarrowtail="3\msx@1B
\mathchardef\leftrightarrows="3\msx@1C
\mathchardef\rightleftarrows="3\msx@1D
\mathchardef\Lsh="3\msx@1E
\mathchardef\Rsh="3\msx@1F
\mathchardef\rightsquigarrow="3\msx@20
\mathchardef\leftrightsquigarrow="3\msx@21
\mathchardef\looparrowleft="3\msx@22
\mathchardef\looparrowright="3\msx@23
\mathchardef\circeq="3\msx@24
\mathchardef\succsim="3\msx@25
\mathchardef\gtrsim="3\msx@26
\mathchardef\gtrapprox="3\msx@27
\mathchardef\multimap="3\msx@28
\mathchardef\therefore="3\msx@29
\mathchardef\because="3\msx@2A
\mathchardef\doteqdot="3\msx@2B

\mathchardef\triangleq="3\msx@2C
\mathchardef\precsim="3\msx@2D
\mathchardef\lesssim="3\msx@2E
\mathchardef\lessapprox="3\msx@2F
\mathchardef\eqslantless="3\msx@30
\mathchardef\eqslantgtr="3\msx@31
\mathchardef\curlyeqprec="3\msx@32
\mathchardef\curlyeqsucc="3\msx@33
\mathchardef\preccurlyeq="3\msx@34
\mathchardef\leqq="3\msx@35
\mathchardef\leqslant="3\msx@36
\mathchardef\lessgtr="3\msx@37
\mathchardef\backprime="0\msx@38
\mathchardef\risingdotseq="3\msx@3A
\mathchardef\fallingdotseq="3\msx@3B
\mathchardef\succcurlyeq="3\msx@3C
\mathchardef\geqq="3\msx@3D
\mathchardef\geqslant="3\msx@3E
\mathchardef\gtrless="3\msx@3F
\mathchardef\sqsubset="3\msx@40
\mathchardef\sqsupset="3\msx@41
\mathchardef\vartriangleright="3\msx@42
\mathchardef\vartriangleleft ="3\msx@43
\mathchardef\trianglerighteq="3\msx@44
\mathchardef\trianglelefteq="3\msx@45
\mathchardef\bigstar="0\msx@46
\mathchardef\between="3\msx@47
\mathchardef\blacktriangledown="0\msx@48
\mathchardef\blacktriangleright="3\msx@49
\mathchardef\blacktriangleleft="3\msx@4A
\mathchardef\vartriangle="3\msx@4D
\mathchardef\blacktriangle="0\msx@4E
\mathchardef\triangledown="0\msx@4F
\mathchardef\eqcirc="3\msx@50
\mathchardef\lesseqgtr="3\msx@51
\mathchardef\gtreqless="3\msx@52
\mathchardef\lesseqqgtr="3\msx@53
\mathchardef\gtreqqless="3\msx@54
\mathchardef\Rrightarrow="3\msx@56
\mathchardef\Lleftarrow="3\msx@57
\mathchardef\veebar="2\msx@59
\mathchardef\barwedge="2\msx@5A
\mathchardef\doublebarwedge="2\msx@5B
\mathchardef\angle="0\msx@5C
\mathchardef\measuredangle="0\msx@5D
\mathchardef\sphericalangle="0\msx@5E
\mathchardef\varpropto="3\msx@5F
\mathchardef\smallsmile="3\msx@60
\mathchardef\smallfrown="3\msx@61
\mathchardef\Subset="3\msx@62
\mathchardef\Supset="3\msx@63
\mathchardef\Cup="2\msx@64

\mathchardef\Cap="2\msx@65

\mathchardef\curlywedge="2\msx@66
\mathchardef\curlyvee="2\msx@67
\mathchardef\leftthreetimes="2\msx@68
\mathchardef\rightthreetimes="2\msx@69
\mathchardef\subseteqq="3\msx@6A
\mathchardef\supseteqq="3\msx@6B
\mathchardef\bumpeq="3\msx@6C
\mathchardef\Bumpeq="3\msx@6D
\mathchardef\lll="3\msx@6E

\mathchardef\ggg="3\msx@6F

\mathchardef\circledS="0\msx@73
\mathchardef\pitchfork="3\msx@74
\mathchardef\dotplus="2\msx@75
\mathchardef\backsim="3\msx@76
\mathchardef\backsimeq="3\msx@77
\mathchardef\complement="0\msx@7B
\mathchardef\intercal="2\msx@7C
\mathchardef\circledcirc="2\msx@7D
\mathchardef\circledast="2\msx@7E
\mathchardef\circleddash="2\msx@7F
\def\ulcorner{\delimiter"4\msx@70\msx@70 }
\def\urcorner{\delimiter"5\msx@71\msx@71 }
\def\llcorner{\delimiter"4\msx@78\msx@78 }
\def\lrcorner{\delimiter"5\msx@79\msx@79 }
\def\yen{{\mathhexbox@\msx@55 }}
\def\checkmark{{\mathhexbox@\msx@58 }}
\def\circledR{{\mathhexbox@\msx@72 }}
\def\maltese{{\mathhexbox@\msx@7A }}
\mathchardef\lvertneqq="3\msy@00
\mathchardef\gvertneqq="3\msy@01
\mathchardef\nleq="3\msy@02
\mathchardef\ngeq="3\msy@03
\mathchardef\nless="3\msy@04
\mathchardef\ngtr="3\msy@05
\mathchardef\nprec="3\msy@06
\mathchardef\nsucc="3\msy@07
\mathchardef\lneqq="3\msy@08
\mathchardef\gneqq="3\msy@09
\mathchardef\nleqslant="3\msy@0A
\mathchardef\ngeqslant="3\msy@0B
\mathchardef\lneq="3\msy@0C
\mathchardef\gneq="3\msy@0D
\mathchardef\npreceq="3\msy@0E
\mathchardef\nsucceq="3\msy@0F
\mathchardef\precnsim="3\msy@10
\mathchardef\succnsim="3\msy@11
\mathchardef\lnsim="3\msy@12
\mathchardef\gnsim="3\msy@13
\mathchardef\nleqq="3\msy@14
\mathchardef\ngeqq="3\msy@15
\mathchardef\precneqq="3\msy@16
\mathchardef\succneqq="3\msy@17
\mathchardef\precnapprox="3\msy@18
\mathchardef\succnapprox="3\msy@19
\mathchardef\lnapprox="3\msy@1A
\mathchardef\gnapprox="3\msy@1B
\mathchardef\nsim="3\msy@1C
\mathchardef\napprox="3\msy@1D
\mathchardef\varsubsetneq="3\msy@20
\mathchardef\varsupsetneq="3\msy@21
\mathchardef\nsubseteqq="3\msy@22
\mathchardef\nsupseteqq="3\msy@23
\mathchardef\subsetneqq="3\msy@24
\mathchardef\supsetneqq="3\msy@25
\mathchardef\varsubsetneqq="3\msy@26
\mathchardef\varsupsetneqq="3\msy@27
\mathchardef\subsetneq="3\msy@28
\mathchardef\supsetneq="3\msy@29
\mathchardef\nsubseteq="3\msy@2A
\mathchardef\nsupseteq="3\msy@2B
\mathchardef\nparallel="3\msy@2C
\mathchardef\nmid="3\msy@2D
\mathchardef\nshortmid="3\msy@2E
\mathchardef\nshortparallel="3\msy@2F
\mathchardef\nvdash="3\msy@30
\mathchardef\nVdash="3\msy@31
\mathchardef\nvDash="3\msy@32
\mathchardef\nVDash="3\msy@33
\mathchardef\ntrianglerighteq="3\msy@34
\mathchardef\ntrianglelefteq="3\msy@35
\mathchardef\ntriangleleft="3\msy@36
\mathchardef\ntriangleright="3\msy@37
\mathchardef\nleftarrow="3\msy@38
\mathchardef\nrightarrow="3\msy@39
\mathchardef\nLeftarrow="3\msy@3A
\mathchardef\nRightarrow="3\msy@3B
\mathchardef\nLeftrightarrow="3\msy@3C
\mathchardef\nleftrightarrow="3\msy@3D
\mathchardef\divideontimes="2\msy@3E
\mathchardef\varnothing="0\msy@3F
\mathchardef\nexists="0\msy@40
\mathchardef\mho="0\msy@66
\mathchardef\thorn="0\msy@67
\mathchardef\beth="0\msy@69
\mathchardef\gimel="0\msy@6A
\mathchardef\daleth="0\msy@6B
\mathchardef\lessdot="3\msy@6C
\mathchardef\gtrdot="3\msy@6D
\mathchardef\ltimes="2\msy@6E
\mathchardef\rtimes="2\msy@6F
\mathchardef\shortmid="3\msy@70
\mathchardef\shortparallel="3\msy@71
\mathchardef\smallsetminus="2\msy@72
\mathchardef\thicksim="3\msy@73
\mathchardef\thickapprox="3\msy@74
\mathchardef\approxeq="3\msy@75
\mathchardef\succapprox="3\msy@76
\mathchardef\precapprox="3\msy@77
\mathchardef\curvearrowleft="3\msy@78
\mathchardef\curvearrowright="3\msy@79
\mathchardef\digamma="0\msy@7A
\mathchardef\varkappa="0\msy@7B
\mathchardef\hslash="0\msy@7D
\mathchardef\hbar="0\msy@7E
\mathchardef\backepsilon="3\msy@7F
\def\Bbb{\relaxnext@\ifmmode\let\next\Bbb@\else
 \def\next{\Err@{Use \string\Bbb\space only in math mode}}\fi\next}
\def\Bbb@#1{{\Bbb@@{#1}}}
\def\Bbb@@#1{\noaccents@\fam\msyfam#1}
\catcode`\@=12

\font\tenmsy=msbm10
\font\sevenmsy=msbm7
\font\fivemsy=msbm5
\font\tenmsx=msam10
\font\sevenmsx=msam7
\font\fivemsx=msam5
\newfam\msyfam

\textfont\msyfam=\tenmsy
\scriptfont\msyfam=\sevenmsy
\scriptscriptfont\msyfam=\fivemsy
\newfam\msxfam

\textfont\msxfam=\tenmsx
\scriptfont\msxfam=\sevenmsx
\scriptscriptfont\msxfam=\fivemsx

%
 
\magnification 1200
\baselineskip=24pt plus 6pt
\def\today{\ifcase\month\or January\or February\or March\or April\or
May\or June\or July\or August\or September\or October\or November\or
December \fi\space\number \day, \number\year}

\define\egg{\Relbar\kern-.30em\Relbar\kern-.30em\Relbar}

\define\p1{P^1_{\delta, \lambda}}

\define\pbf{\par\bigpagebreak\flushpar}
\define\cof{\hbox{\rm cof}}


\def\b{\beta}
\def\a{\alpha}
\def\l{\lambda}

\def\ov{\overline}

\def\rest{\mathord\restriction}



\def\k{\kappa}
\magnification=1200
\baselineskip=24pt plus 6pt
\define\cU{\Cal U}



\def\b{\beta}
\def\a{\alpha}
\def\l{\lambda}

\def\ov{\overline}

\def\rest{\mathord\restriction}


\def\k{\kappa}
\def\d{\delta}
\def\no{\noindent}
\def\g{\gamma}
\def\la{\langle}
\def\ra{\rangle}
\def\rest{\hskip 4pt\grave{\!}\!\!\hskip 1pt | \hskip 1 pt}

\def\B{{\cal B}}

\def\u{{\cal U}}

\def\gch{{\hbox{\rm GCH}}}

\def\ath{\alpha^{\hbox{\rm th}}}
\def\bth{\beta^{\hbox{\rm th}}}
\def\pkl{P_\kappa(\lambda)}
\def\id{{\hbox{\rm id}}}

\def\f{\varphi}

\def\pic{\underset i \in C \to{\prod} G_i}
\def\pic0{\underset i \in C_0 \to{\prod} G_i}

\def\sqr#1#2{{\vcenter{\vbox{\hrule height.#2pt
       \hbox{\vrule width.#2pt height#1pt \kern#1pt
         \vrule width.#2pt}
        \hrule height.#2pt}}}}
\def\square{\mathchoice\sqr34\sqr34\sqr{2.1}3\sqr{1.5}3}

\def\hb{\hfil\break}
\def\mag{ \subset_{{}_{{}_{\!\!\!\!\!\sim}}} }

\voffset=-.15in\hoffset=0.00in\vsize=7.65in

\no\S0 Introduction and Preliminaries

Since Solovay defined the notion of supercompact cardinal
in the late 1960s (see [SRK]), ascertaining the nature
of the relationship between supercompact and strongly
compact cardinals has been a prime focus of large
cardinal set theorists. At first, Solovay believed that
every strongly compact cardinal must also be
supercompact. This was refuted by his student Menas in
the early 1970s, who showed in his thesis [Me] that if
$\k$ is the least measurable limit of strongly compact
cardinals, then $\k$ is strongly compact but not
$2^\k$ supercompact. (That this result is best possible
was established about twenty years later by
Shelah and the author.
See [AS$\infty$b] for more details.) Menas further
showed in his thesis [Me] from a measurable limit of
supercompact cardinals that it was consistent for the
least strongly compact cardinal not to be the least
supercompact cardinal. In addition, in unpublished work
that used Menas' ideas, Jacques Stern showed, from
hypotheses on the order of a supercompact limit of
supercompact cardinals, that it was consistent for
the first two strongly compact cardinals not to be
supercompact.

Shortly after Menas' work, Magidor in his celebrated
paper [Ma76] established the fundamental results concerning
the nature of the least strongly compact cardinal,
showing that it was consistent, relative to the consistency
of a strongly compact cardinal, for the least strongly
compact cardinal to be the least measurable cardinal
(in which case, it is not the least supercompact cardinal),
but that it was also consistent, relative to the consistency
of a supercompact cardinal, for the least strongly compact
cardinal to be the least supercompact cardinal. In
generalizations of the above work, Kimchi and Magidor
[KiM] later showed, relative to a class of supercompact
cardinals, that it was consistent for the classes of
supercompact and strongly compact cardinals to coincide,
except at measurable limit points, and for $n \in \omega$,
relative to the consistency of $n$ supercompact cardinals,
it was consistent for the first $n$ measurable cardinals
to be the first $n$ strongly compact cardinals. Further
generalizations of these results can be found in
[A80], [A81], [A95], [A$\infty$a], [A$\infty$b],
[AG], and [AS$\infty$a].

The purpose of this paper is to show that the ideas of
[Me] can be used to force over a model given by
[A$\infty$a] to produce models in which, roughly speaking,
the class of compact cardinals, where a
compact cardinal will be taken as one which is either
strongly compact or supercompact, can have virtually
arbitrary structure. Specifically, we prove the
following two theorems.

\proclaim{Theorem 1}
Let $V \models ``$ZFC + $\Omega$ is the least
inaccessible limit of measurable limits of
supercompact cardinals + $f : \Omega \to 2$ is
a function''. There is then a partial ordering
$P \in V$ so that for $\ov V = V^P$,
$\ov V_\Omega \models ``$ZFC +
There is a proper class of compact cardinals +
If $f(\a) = 0$, then the $\a^{\hbox{\rm th}}$
compact cardinal isn't supercompact +
If $f(\a) = 1$, then the $\a^{\hbox{\rm th}}$
compact cardinal is supercompact''.
\endproclaim

\proclaim{Theorem 2}
Let $V \models ``$ZFC + $\k$ is a supercompact limit
of supercompact cardinals +
$f : \k \to 2$ is a function''.
There is then a partial ordering $P \in V$ so that
$V^P \models ``$ZFC + If $\a$ is not in $V$ a measurable limit
of measurable limits of supercompact cardinals and
$f(\a) = 0$, then the $\a^{\hbox{\rm th}}$
compact cardinal isn't supercompact +
If $\a$ is not in $V$ a measurable limit of measurable
limits of supercompact cardinals and
$f(\a) = 1$, then the $\a^{\hbox{\rm th}}$
compact cardinal is supercompact''.
Further, for any $\a < \k$ which was in $V$
a regular limit of measurable limits of
supercompact cardinals,
$V \models ``\a$ is measurable'' iff
$V^P \models ``\a$ is measurable'', and
every cardinal $\a \le \k$ which was in
$V$ a supercompact limit of supercompact cardinals
remains in $V^P$ a supercompact cardinal.
\endproclaim

We note that in Theorem 2 above, we will have no
control over measurable limits of compact cardinals
in the generic extension. This is since by Menas'
aforementioned result, many of these cardinals $\k$
are provably not $2^\k$ supercompact.

Theorems 1 and 2 have a number of interesting corollaries.
We list a few of these now.

\no 1$.$ In Theorem 1, if $f$ is constantly $0$, then
$\ov V_\Omega \models ``$There is a proper class of
strongly compact cardinals, and no strongly compact
cardinal is supercompact''.

\no 2$.$ In Theorem 1, if $f(\a) = 0$ for even and
limit ordinals, and $f(\a) = 1$ otherwise, then
$\ov V_\Omega \models ``$The compact cardinals alternate
in the pattern non-supercompact, supercompact,
non-supercompact, supercompact, etc$.$, with the
$\a^{\hbox{\rm th}}$
compact cardinal for $\a$ a limit ordinal always
being non-supercompact''.

\no 3$.$ In Theorem 2, if $f$ is as in the last
corollary above, then
$V^P \models ``$The least measurable limit of
supercompact cardinals is the same as the least
measurable limit of non-supercompact strongly
compact cardinals''.

\no Although this corollary easily follows from
Theorem 2, all we will need to prove it is a
model with a measurable limit of measurable limits
of supercompact cardinals.

\no 4$.$ In Theorem 2, if
$f$ is constantly $0$, then
$V^P \models ``$The least supercompact cardinal is
a limit of strongly compact cardinals''.

\no We will indicate (with some details missing)
following the proof of Theorem 2 how Corollary 4
is proven and how Corollary 3 is proven using the
weaker hypotheses mentioned above.

The structure of this paper is as follows.
Section 0 contains our Introduction and Preliminaries.
Section 1 contains the proof of Theorem 1.
Section 2 contains the proof of Theorem 2.
Section 3 contains a discussion of the proofs of
Corollaries 3 and 4 and some concluding remarks.

We digress now to give some preliminary information.
Essentially, our notation and terminology are standard, and
when this is not the
case, this will be clearly noted.
For $\a < \b$ ordinals, $[\a, \b], [\a, \b), (\a, \b]$, and $(\a,
\b)$ are as
in standard interval notation.

When forcing, $q \ge p$ will mean that $q$ is stronger than $p$,
and for $\f$ a formula in the forcing language with
respect to our partial ordering $P$ and $p \in P$, $p \Vert \f$
will mean that
$p$ decides $\f$.
For $G$ $V$-generic over $P$, we will use both $V[G]$ and $V^{P}$
to indicate the universe obtained by forcing with $P$.
If $x \in V[G]$, then $\dot x$ will be a term in $V$ for $x$.
We may, from time to time, confuse terms with the sets they denote
and write $x$
when we actually mean $\dot x$, especially
when $x$ is some variant of the generic set $G$, or $x$ is
in the ground model $V$.

If $\k$ is a cardinal and $P$ is
a partial ordering, $P$ is $\k$-closed if given a sequence
$\langle p_\a: \a < \k \rangle$ of elements of $P$ so that
$\beta < \gamma < \k$ implies $p_\beta \le p_\gamma$ (an increasing
chain of length
$\k$), then there is some $p \in P$ (an upper bound to this chain)
so that $p_\a \le p$ for all $\a < \k$.
$P$ is $<\k$-closed if $P$ is $\delta$-closed for all cardinals
$\d < \k$.
$P$ is $\k$-directed closed if for every cardinal $\delta < \k$
and every directed
set $\langle p_\a : \a < \delta \rangle $ of elements of $P$
(where $\langle p_\alpha : \alpha < \delta \rangle$ is directed if
for every two distinct elements $p_\rho, p_\nu \in
\langle p_\alpha : \alpha < \delta \rangle$, $p_\rho$ and
$p_\nu$ have a common upper bound) there is an
upper bound $p \in P$. $P$ is $\k$-strategically closed if in the
two person game in which the players construct an increasing
sequence
$\langle p_\a: \a \le\k\rangle$, where player I plays odd
stages and player
II plays even and limit stages, then player II has a strategy which
ensures the game can always be continued.
Note that if $P$ is $\k$-strategically closed and
$f : \k \to V$ is a function in $V^P$, then $f \in V$.
$ P$ is $< \k$-strategically closed if $P$ is $\delta$-strategically
closed for all cardinals $\delta < \k$.
$P$ is $\prec   \k$-strategically closed if in the two
person game in which the players construct an increasing
sequence $\langle p_\alpha : \alpha < \k \rangle$, where
player I plays odd stages and player II plays even and limit
stages, then player II has a strategy which ensures the game
can always be continued.
Note that trivially, if $P$ is $<\k$-closed, then $P$ is
$< \k$-strategically
closed and $\prec \k  $-strategically closed. The converse of
both of these facts is false.

We mention that we are assuming complete familiarity with the
notions of
measurability, strong compactness, and supercompactness.
Interested readers may consult [SRK], [Ka], or [KaM] for
further details.
We note first that all elementary embeddings witnessing the $\l$
supercompactness  of $\k$ will come from some
fine, $\k$-complete, normal
ultrafilter $\cU$ over $P_\k (\l) = \{ x \subseteq \l: | x| < \k \}
$, and all elementary embeddings witnessing the $\l$
strong compactness of $\k$ will come from some
fine, $\k$-complete ultrafilter ${\cal U}$ over
$P_\k(\l)$.

We note also the following properties, which
will be used throughout the course of the paper.

\no 1. (Menas [Me]) If $\k$ is the $\ath$
measurable limit of strongly compact or supercompact
cardinals and $\a < \k$, then $\k
$ is strongly compact but isn't $2^\k$ supercompact.
A proof of this fact for the $\ath$ measurable limit of
strongly compact
cardinals will be given during the proof of Lemma 4.
The proof for the $\ath$ measurable limit of supercompact
cardinals is the same.

\no 2. ([SRK]) If $\d < \k \le \l$ are regular cardinals
and $\k$ is strongly compact, then every stationary subset
$S \subseteq \l$ of ordinals of cofinality $\d$ reflects, i.e.,
for some ordinal $\a < \l$,
$S \cap \a$ is stationary at its supremum.

\no 3. (Magidor [Ma71]) If $\k < \l$ are so that $\k$
is $< \l$ supercompact and $\l$ is supercompact, then
$\k$ is supercompact.

\no 4. (DiPrisco [DH]) If $\k < \l$ are so that $\k$
is $< \l$ strongly compact and $\l$ is strongly compact,
then $\k$ is strongly compact.

Let $\g < \k$ be so that $\g$ and $\k$ are
regular cardinals. We now describe and state the properties
of the standard notion of forcing $P_{\g, \k}$
for adding a non-reflecting stationary set of ordinals of
cofinality $\g$ to $\k$. Specifically,
$P_{\g, \k} = \{p$ : For some $\a < \k$,
$p : \a \to \{0,1\}$ is a characteristic function of
$S_p$, a subset of $\a$ not stationary at its
supremum nor having any initial segment which is
stationary at its supremum, so that $\b \in S_p$ implies
$\b > \g$ and cof$(\b) = \g \}$, ordered by $q \ge p$
iff $q \supseteq p$ and
$S_p = S_q \cap \sup(S_p)$, i.e.,
$S_q$ is an end extension of $S_p$. It is well-known
that for $G$ $V$-generic over $P_{\g, \k}$
(see [Bu] or [KiM]),
in $V[G]$, a non-reflecting stationary set
$S = S[G] = \cup \{S_p : p \in G \} \subseteq \k$
of ordinals of cofinality $\g$ has been introduced, and
since $P_{\g, \k}$ is
$\prec \k$-strategically closed,
in $V[G]$, the
bounded subsets of $\k$ are the same as those in $V$.
It is also virtually immediate that $P_{\g, \k}$ is
$\g$-directed closed.

It is clear from the definition of $P_{\g, \k}$
that assuming GCH holds in our ground model $V$,
$|P_{\g, \k}| = \k$. Thus, the strategic closure
properties of $P_{\g, \k}$ mentioned in the above
paragraph imply
$V^{P_{\g, \k}} \models \gch$.
Also, if $\la \k_\a : \a < \l \ra$
is a strictly increasing sequence of regular cardinals
and $\la \g_\a : \a < \l \ra$ is a sequence of regular
cardinals (not necessarily distinct) so that
$\g_\a < \k_\a$ for all $\a < \l$, then if
$P = \la \la P_\a, \dot Q_\a \ra : \a < \l \ra$
is the Easton support iteration where
$P_0 = \{\emptyset\}$ and
$\force_{P_\a} ``\dot Q_\a = \dot P_{\g_\a, \k_\a}$'',
then since Easton support iterations of strategically
closed partial orderings retain the appropriate amount of
strategic closure, the standard arguments in
combination with the above mentioned cardinality and
strategic closure properties imply
$V^P \models \gch$.
Further, if
$R^* = \la R_\a : \a < \d \ra$ is a sequence of
partial orderings where each $R_\a$ is an
iteration as described in the preceding sentence
and $R^*$ is so that for $\b_\a$
the sup of the cardinals in the domain of $R_\a$,
$0 \le \a_0 < \a_1 < \d$ implies
$\b_{\a_0} < \b_{\a_1}$, then for $R$ the Easton support product
$\underset \a < \d \to{\prod} R_\a$,
it is once more the case that
$V^R \models \gch$.

\no\S1 The Proof of Theorem 1

We turn now to the proof of Theorem 1. Recall we are assuming
$V \models ``$ZFC + $\Omega$ is the least inaccessible limit
of measurable limits of supercompact cardinals +
$f : \Omega \to 2$ is a function''. By the results of
[A$\infty$a], we also assume without loss of generality that
$V \models ``$The supercompact and strongly compact cardinals
coincide except at measurable limit points + Every
supercompact cardinal is Laver indestructible [L]''.

Before defining the partial ordering $P$ used in the
proof of Theorem 1, we fix first some notation to be used
throughout the duration of the proof of Theorem 1.
For $\a < \Omega$, let $\d_\a$ be the $\ath$ measurable
limit of supercompact cardinals
(which, since $\a < \d_\a$, means by Menas' result
stated above that $\d_\a$ isn't
$2^{\d_\a}$ supercompact), and let
$\la \k^\a_\b : \b < \d_\a \ra$
be an increasing sequence of supercompact cardinals whose limit is
$\d_\a$ so that
$\k^\a_0 > \underset \b < \a \to{\cup} \d_\b$
and so that all supercompact cardinals in the interval
$(\underset \b < \a \to{\cup} \d_\b, \d_\a)$
(which in this instance is the same as all supercompact
cardinals in the interval
$[\underset \b < \a \to{\cup} \d_\b, \d_\a)$,
since $\underset \b < \a \to{\cup} \d_\b$
isn't supercompact, being below the least measurable
limit of measurable limits of supercompact cardinals)
are elements of
$\la \k^\a_\b : \b < \d_\a \ra$.

We define now the partial ordering $P$ used in the
proof of Theorem 1.
If $f(\a) = 0$, $P_\a$ is the Easton support iteration
$\la Q_\b \ast \dot R_\b : \b < \d_\a \ra$,
where $Q_0 = \{\emptyset\}$ and
$\force_{Q_\b} ``\dot R_\b$ adds a non-reflecting
stationary set of ordinals of cofinality
${(\underset \g < \a \to{\cup} \d_\g)}^+$
to $\k^\a_\b$''.
If $f(\a) = 1$, $P_\a$ is the partial ordering for
adding a non-reflecting stationary set of ordinals of
cofinality $\k^\a_0$ to $\d_\a$.
The partial ordering $P$ used in the proof of Theorem 1
is then defined as the Easton support product
$\underset \a < \Omega \to{\prod} P_\a$.

The intuition behind the definition of $P$ is quite
simple. If the $\ath$ compact cardinal
in our final model $\ov V_\Omega$ is to be
non-supercompact, then we start with the
$\ath$ measurable limit of supercompact
cardinals $\d_\a$, a cardinal which is provably
strongly compact but not supercompact, and destroy
all supercompact cardinals below $\d_\a$ but beyond
$\underset \b < \a \to{\cup} \d_\b$.
Since we start with a model in which the strongly
compact and supercompact cardinals coincide except at
measurable limit points, we will have after forcing that
the $\ath$ compact cardinal isn't
supercompact and has no compact cardinals below it
except for those explicitly preserved by the forcing.
If, however, the $\ath$ compact cardinal in
$\ov V_\Omega$ is to be supercompact, then we
destroy the strong compactness of $\d_\a$ by a forcing
which will preserve the supercompactness of
$\k^\a_0$ and the strong compactness and supercompactness
of those cardinals below $\k^\a_0$ which are to become in
$\ov V_\Omega$ the compact cardinals below $\k^\a_0$,
yet will destroy all compact cardinals in the interval
$(\k^\a_0, \d_\a]$.

\proclaim{Lemma 1}
If $f(\a) = 1$,
$V^P \models ``\k^\a_0$ is supercompact''.
\endproclaim

\demo{Proof of Lemma 1}
Write now and for the rest of the proofs of Theorems 1 and 2
$P$ as
$P^\a \times P_\a \times P_{< \a}$,
where $P_{< \a}$ and $P^\a$ are Easton support products,
$P_{< \a} = \underset \b < \a \to{\prod} P_\b$,
and
$P^\a = \underset \b \in (\a, \Omega) \to{\prod} P_\a$.
By the definition of $P_\b$ for
$\b \in [\a, \Omega)$, $P^\a \times P_\a$ is
$\k^\a_0$-directed closed. Therefore, since
$V \models ``\k^\a_0$ is Laver indestructible'',
$V^{P^\a \times P_\a} \models ``\k^\a_0$ is
supercompact''. Also, since $\Omega$ is the least
inaccessible limit of measurable limits of supercompact
cardinals, $|P_{< \a}| < \k^\a_0$.
Thus, by the L\'evy-Solovay results [LS],
$V^{P^\a \times P_\a \times P_{< \a}} = V^P \models
``\k^\a_0$ is supercompact''.
This proves Lemma 1.
\pbf
\hfill $\square$ Lemma 1

Our next goal will be to show that if $f(\a) = 0$, then
$V^P \models ``\d_\a$ is a non-supercompact
strongly compact cardinal''.
This will be done using ideas of Menas found in [Me].
Before doing this, however, we will prove two technical lemmas.
The first is a lemma
of Menas about the existence of certain kinds of
strongly compact ultrafilters over $P_\k(\l)$ when
$\k$ is a measurable limit of strongly compact cardinals.
The second shows that if $\k$ is strongly compact
in $V$ and
$Q$ is a partial ordering so that $V$ and $V^Q$
contain the same bounded subsets of $\k$, then any
strongly compact cardinal in $V^Q$ below $\k$ is also
strongly compact in $V$.

\proclaim{Lemma 2 (Proposition 2.31 of [Me])}
Let $\k < \l$ be cardinals with $\k$ a measurable
limit of strongly compact cardinals.
Let $f' : \k \to \k$ be defined by
$f'(\a) =$ The least strongly compact cardinal above $\a$.
There is then a strongly compact ultrafilter $\u$
over $P_\k(\l)$ so that for
$j_{\u} : V \to M_{\u}$
the associated strongly compact elementary embedding and $g$
the function representing $\k$ in $M_{\u}$,
$\{ p \in \pkl : f'(g(p)) > |p| \} \in \u$.
\endproclaim

\demo{Proof of Lemma 2}
Let $\mu$ be a normal measure over $\k$. Define
$f'' : \k \to \k$ by
$f''(\a) =$ The sup of all strongly compact cardinals
below $\a$.
It is clear $\{ \a : f''(\a) \le \a \} \in \mu$. If
$\{ \a : f''(\a) < \a \} \in \mu$, then by the
normality of $\mu$, $\{ \a : f''(\a) = \a_0 \} \in \mu$
for some $\a_0 < \k$. This, however, contradicts the fact
that $\k$ is a limit of strongly compact cardinals, so
$A = \{ \a < \k : \a$ is a limit of strongly compact
cardinals$\} \in \mu$. This means that for $\a < \b$
in $A$, $\a, \b$ arbitrary, $f'(\a) < \b$.

For every $\a \in A$, let $\mu_\a$ be a
strongly compact ultrafilter over
$P_{f'(\a)}(\l)$. Let
$\u$ be defined by
$X \in \u$ iff $X \subseteq P_\k(\l)$ and
$\{ \a < \k : X \cap
P_{f'(\a)}(\l) \in \mu_\a \} \in \mu$.
It is easily checked that $\u$ is a strongly
compact ultrafilter over $\pkl$. We show that
$\u$ has the desired property.

For every $\a \in A$, let $B_\a = \{ p \in
P_{f'(\a)}(\l) : |p| \in (\a, f'(\a)) \}$.
By the fineness of $\mu_\a$, $B_\a \in \mu_\a$, so
$B = \underset \a \in A \to{\cup} B_\a \in \u$.
Also, by the choice of $A$, for every $p \in B$,
there is a unique $\a \in A$ so that $p \in B_\a$.
This means the function
$g(p) =$ The unique $\a \in A$ so that $p \in B_\a$
is well-defined for $p \in B$. It is again clear by
the first sentence of this paragraph that for every
$p \in B$, $f'(g(p)) > |p|$, i.e.,
$\{ p \in \pkl : f'(g(p)) > |p| \} \in \u$.
Thus, the proof of Lemma 2 will be complete once we have shown
${[g]}_{\u} = \k$.

To show this last fact, let $h$ be so that
$\{ p \in \pkl : h(p) < g(p) \} \in \u$.
This means by the definition of $\u$ and the fact
$B \in \u$ that we may assume for some $C \subseteq A$,
$C \in \mu$, for every $\a \in C$, $B'_\a =
\{ p \in B_\a : h(p) < g(p) \} \in \mu_\a$.
Let $\a \in C$ be arbitrary. Since for
$p \in B'_\a \subseteq B_\a$, $|p| \in (\a, f'(\a))$
and $g(p) = \a$, for $p \in B'_\a$,
$h(p) < g(p) = \a < f'(\a)$. Thus, for some
$B''_\a \subseteq B'_\a$, $B''_\a \in \mu_\a$,
the additivity of $\mu_\a$ implies the existence of a
$\b_\a < \a$ so that for every $p \in B''_\a$,
$h(p) = \b_\a$. If we now define
$h' : C \to \k$ by $h'(\a) = \b_\a$, then
$h'(\a) < \a$ for all $\a \in C$. Thus, by the
normality of $\mu$, for some $D \subseteq C$,
$D \in \mu$ and some fixed $\b < \k$, $\a \in D$
implies $h(p) = \b$ for every $p \in B''_\a$.
This means that $\{ p \in \pkl : h(p) = \b \}
\in \u$. Since for any fixed $\g < \k$,
$\{ p \in \pkl : g(p) > \g \} \in \u$, we can now infer that
${[g]}_{\u} = \k$. This proves Lemma 2.
\pbf
\hfill $\square$ Lemma 2

We remark that the referee has pointed out an alternative
proof of Lemma 2 is possible using elementary embeddings.
An outline of the argument is as follows, where we adopt
the notation of Lemma 2.
Let $j_\mu : V \to M_\mu$ be the ultrapower embedding given by
$\mu$.
There is then $k : M_\mu \to N$ witnessing that
$\k^* = j_\mu(f')(\k)$ is $j_\mu(\l)$ strongly compact, so let
$X \in N$ be so that $k''j_\mu(\l) \subseteq X$ and
$N \models ``|X| < k(\k^*)$''. It is easily verifiable that
$k \circ j_\mu$ witnesses the $\l$ strong compactness
of $\k$. If $X$ is chosen so that the $\l$ strong compactness
measure ${\cal U} = \{Z : X \in k \circ j_\mu(Z) \}$
is such that for $j_{\cal U} : V \to M_{\cal U}$
the ultrapower embedding,
$k \circ j_\mu = j_{\cal U}$ and
$M_{\cal U} = N$, then ${\cal U}$ has the desired property.

\proclaim{Lemma 3}
Suppose $V \models ``\k$ is strongly compact'' and $Q$ is a
partial ordering so that $V$ and $V^Q$ contain the
same bounded subsets of $\k$. Then for $\sigma < \k$,
if $V^Q \models ``\sigma$ is strongly compact'',
$V \models ``\sigma$ is strongly compact''.
\endproclaim

\demo{Proof of Lemma 3}
Since $V$ and $V^Q$ contain the same bounded subsets of
$\k$ (meaning $\k$ is a strong limit cardinal in both
$V$ and $V^Q$),
$V \models ``\sigma$ is $< \k$ strongly compact''.
Thus, by the theorem of DiPrisco [DH] mentioned in
the Introduction,
$V \models ``\sigma$ is strongly compact''.
This proves Lemma 3.
\pbf
\hfill $\square$ Lemma 3

\proclaim{Lemma 4}
If $f(\a) = 0$,
$V^P \models ``\d_\a$ is strongly compact''.
\endproclaim

\demo{Proof of Lemma 4}
The definition of $P_\b$ for $\b \in (\a, \Omega)$
implies each $P_\b$ for $\b \in (\a, \Omega)$ must be
at least $\d^+_\a$-directed closed. Thus, $P^\a$ is at
least $\d^+_\a$-directed closed. Therefore, since
$V \models ``$All supercompact cardinals are Laver
indestructible'',
$V^{P^\a} \models ``\d_\a$ is a measurable limit of
supercompact cardinals'', i.e.,
$V^{P^\a} \models ``\d_\a$ is strongly compact''.

Call $V^{P^\a}$ $V^0$ and $\d_\a$ $\d$. We show now that
${(V^0)}^{P_\a} \models ``\d$ is strongly compact''.
The proof is essentially the same as the proof of
Theorem 2.27 of [Me]. Let $\g \ge \d$ be arbitrary, and let
$\l =
2^{[\g]^{< \d}}
$. Let $\u$ be a strongly compact ultrafilter over $
P_\d(\l)
$ having the property of Lemma 2, and let
$j : V^0 \to M$ be the associated strongly compact
elementary embedding.

We begin by noting that $M \models ``\d$ is not
measurable''. To see this, we remark first that
$V^0 \models ``\d$ is the $\ath$ measurable limit of
strongly compact cardinals''. To prove this last
fact about $V^0$, observe that
$V \models ``\d$ is the $\ath$ measurable limit of
strongly compact cardinals'', and as already noted,
$P^\a$ is $\d^+$-directed closed.
This means $V$ and $V^0$ contain the same bounded
subsets of $\d$ and $V^0 \models ``\d$ is
measurable''. Thus, by Lemma 3, any strongly compact
cardinal in $V^0$ below $\d$ is already strongly
compact in $V$, so $V^0 \models ``\d$ is the
$\ath$ measurable limit of strongly compact
cardinals''.

The rest of the argument
that $M \models ``\d$ is not measurable''
parallels the argument given
in Lemma 12 of [AS$\infty$b] (which is different from
the argument Menas gives in Theorem 2.22 of
[Me]). If
$M \models ``\d$ is measurable'', then
since $\a < \d$ and $j \rest \d = \id$,
$M \models ``\d$ is the $j(\ath) = \ath$ measurable
limit of strongly compact cardinals''. This,
of course, contradicts that
$j(\d) > \d$ and
$M \models ``j(\d)$ is the $j(\ath) = \ath$
measurable limit of strongly compact cardinals''. Thus,
$M \models ``\d$ is not measurable''. This means that
in $M$, $j(P_\a) = P_\a \ast \dot Q$, where $\d$ is
not in the domain of $\dot Q$. Further, by the
definition of $P_\a$ in both $V$ and $V^0$ and
the property of $\u$ given by Lemma 2, in $M$,
the least cardinal $\sigma$ in the domain of $\dot Q$
is so that $\sigma > |
{[\id]}_{\u}
|$.

Let $G$ be $V^0$-generic over $P_\a$, and let
$H$ be $V^0[G]$-generic over $Q$. By the above
factorization property of $j(P_\a)$ in $M$,
$j : V^0 \to M$ extends in the usual way in
$V^0[G \ast H]$ to the elementary embedding
$j^* : V^0[G] \to M[G \ast H]$ given by
$j^*(i_G(\tau)) = i_{G \ast H}(j(\tau))$. $j^*$
can then be used in $V^0[G \ast H]$ to define the set
$\mu$ given by $X \in \mu$ iff $X \subseteq
{(P_\d(\g))}^{V^0[G]}
$ and ${[\id \rest \g]}_{\u} \in j^*(X)$, where
$\id \rest \g : P_\d(\l) \to P_\d(\g)$ is the function
$\id \rest \g(p) = p \cap \g$.
It is easy to check (and is left as an exercise for readers)
that $\mu$ defines, in $V^0[G \ast H]$, a
strongly compact ultrafilter over $
{(P_\d(\g))}^{V^0[G]}
$. We will be done once we have shown $\mu \in V^0[G]$.

To do this, let in $V^0$ $g' : \l \to \dot r$
be a surjection, where $\dot r$ is so that $i_G(\dot r) =
{(2^{[\g]^{< \d}})}^{V^0[G]}$.
(The choice of $\l$ ensures such a surjection exists.)
Let $g$ be a function defined on $P_\d(\l)$ so that
$g(p) = g' \rest p$. Then
$M \models ``j(g)$ is a function from
${[\id]}_{\u}$
into $j(\dot r)$''. This allows us to define a function
$h :
 {[\id]}_{\u} 
\to 2$ in $M[G \ast H]$ by $h(x) = 1$ iff
${[\id \rest \g]}_{\u} \in i_{G \ast H}(j(g)(x))$.
Since the least element $\sigma$ in the domain of
$\dot Q$ is $> |{[\id]}_{\u}|$, and since
by the definition of $P_\a$,
$M[G] \models ``Q$ is $< \sigma$-strategically closed'',
it is the case that
$h \in M[G] \subseteq V[G]$, i.e., $h \in V[G]$.
And, as can be verified, for every $\a < \l$,
$i_G(g'(\a)) \in \mu$ iff for some $q \in G$,
$q \force ``g'(\a) \subseteq
{(P_\d(\g))}^{V^0[G]}$''
and $h(j(\a)) = 1$. This immediately implies
$\mu \in V^0[G]$. Thus,
$V^{P^\a \times P_\a} \models ``\d_\a$ is
strongly compact''. Therefore, since the definition of
$P$ ensures that as in Lemma 1,
$V \models ``|P_{< \a}| < \d_\a$'', the arguments
of [LS] once again tell us
$V^{P^\a \times P_\a \times P_{< \a}} = V^P \models
``\d_\a$ is strongly compact''. This proves Lemma 4.
\pbf
\hfill $\square$ Lemma 4

\proclaim{Lemma 5}
If $f(\a) = 0$,
$\ov V = V^P \models ``\d = \d_\a$ isn't supercompact''.
\endproclaim

\demo{Proof of Lemma 5}
By Lemma 4, for $\l \ge \d$ arbitrary, we can fix
$j : \ov V \to M$ to be an elementary embedding
witnessing the $\l$ strong compactness of $\d$.
Since $\a < \d$ and
$V \models ``\d$ is the $\ath$ measurable limit of
strongly compact cardinals'', for some $\b \le \a$,
$\ov V \models ``\d$ is the $\bth$ measurable
cardinal so that in $V$, $\d$ is a measurable limit of
strongly compact cardinals''. By elementariness and the
facts $\b \le \a < \d$ and $j \rest \d = \id$, if
$M \models ``\d$ is measurable'',
$M \models ``\d$ is the $j(\bth) = \bth$
cardinal so that in $j(V)$, $\d$ is a measurable
limit of strongly compact cardinals''. This, of course,
contradicts that
$M \models ``j(\d) > \d$ is the
$j(\bth) = \bth$ cardinal so that in $j(V)$,
$j(\d)$ is a measurable limit of strongly compact cardinals''.
Thus, $j$ can't be an embedding witnessing the
$2^\d$ supercompactness of $\d$. This proves Lemma 5.
\pbf
\hfill $\square$ Lemma 5

\proclaim{Lemma 6}
$\ov V = V^P \models ``$If $f(\a) = 0$, $\d_\a$
is the $\ath$ strongly compact cardinal, but if
$f(\a) = 1$, $\k^\a_0$ is the $\ath$ strongly
compact cardinal''.
\endproclaim

\demo{Proof of Lemma 6}
Assume Lemma 6 is true for all $\b < \a$. By Lemmas 1 and 4,
the definition of $P_\g$ for any $\g$,
and the fact the theorem of [SRK] mentioned in the
Introduction tells us that if
$\rho$ contains a non-reflecting stationary set of ordinals
of cofinality $\sigma$, then there are no
strongly compact cardinals in the interval $(\sigma, \rho]$,
if $f(\a) = 0$,
$V^{P^\a \times P_\a} \models
``\d_\a$ is strongly compact and there are no strongly compact
cardinals in the interval $[\k^\a_0, \d_\a)$'',
but if $f(\a) = 1$,
$V^{P^\a \times P_\a} \models
``\k^\a_0$ is supercompact and there are no strongly
compact cardinals in the interval $(\k^\a_0, \d_\a]$''. If
$\zeta = \underset \b < \a \to{\cup} \d_\b$,
then by the definition of $P_{< \a}$,
$|P_{< \a}| \le 2^\zeta < \k^\a_0$.
Further, since in $V$, the strongly compact and supercompact
cardinals coincide except at measurable limit points,
the definition of $\zeta$ tells us
$V \models ``$There are no strongly compact cardinals in the
interval $(\zeta, \k^\a_0)$''. By
Lemma 3, since $V$ and
$V^{P^\a \times P_\a}$ have the same bounded subsets of
$\k^\a_0$,
$V$ and $V^{P^\a \times P_\a}$ have the same strongly compact
cardinals $< \k^\a_0$. The arguments of [LS] then yield that
$V^{P^\a \times P_\a \times P_{< \a}} = V^P \models
``$There are no strongly compact cardinals in the interval
$(\zeta, \k^\a_0)$''. This immediately allows us to conclude that
$V^P \models ``$If $f(\a) = 0$, $\d_\a$ is the $\ath$
strongly compact cardinal, but if $f(\a) = 1$,
$\k^\a_0$ is the $\ath$ strongly compact cardinal''.
This proves Lemma 6.
\pbf
\hfill $\square$ Lemma 6

\proclaim{Lemma 7}
$V^P \models ``\Omega$ is inaccessible''.
\endproclaim

\demo{Proof of Lemma 7}
As indicated in the proof of Lemma 4, for any
$\a < \Omega$, $P^\a$ is at least
$\d^+_\a$-directed closed. Further, regardless if
$f(\a) = 0$ or $f(\a) = 1$, by the definition of $P$,
$P_\a$ is $< \k^\a_0$-strategically closed and
$|P_{< \a}| < \k^\a_0$. Thus, since $\Omega$ is
regular in $V$,
$V^{P^\a \times P_\a \times P_{< \a}} = V^P \models
``\cof(\Omega) \ge \k^\a_0$''. As the $\k^\a_0$ are
unbounded in $\Omega$,
$V^P \models ``\Omega$ is regular''. And, by Lemma 6,
$\Omega$ is in $V^P$ a limit of compact cardinals,
meaning
$V^P \models ``\Omega$ is a strong limit cardinal''. Thus,
$V^P \models ``\Omega$ is inaccessible''. This proves Lemma 7.
\pbf
\hfill $\square$ Lemma 7

Lemmas 1-7 complete the proof of Theorem 1.
\pbf
\hfill $\square$ Theorem 1

We remark that it is possible to get sharp bounds in
$V$, $\ov V$, and $\ov V_\Omega$ on the
non-supercompactness of each $\d_\a$ for which
$f(\a) = 0$. We may assume by the methods of
[A$\infty$a] that in the ground model $V$, GCH holds
and each supercompact cardinal $\k$ has been made
indestructible only under forcing with $\k$-directed
closed partial orderings not destroying GCH. This tells
us GCH holds in both $\ov V$ and $\ov V_\Omega$.
It will then be the case by the arguments given in
Lemma 12 of [AS$\infty$b]
(which were used in the fourth paragraph of the
proof of Lemma 4)
that for each $\d < \Omega$ so that
$V \models ``\d$ is a measurable limit of strongly
compact cardinals'',
$V \models ``\d$ isn't $2^\d = \d^+$ supercompact''.
Therefore, since GCH holds in both
$\ov V$ and $\ov V_\Omega$, as observed in the proof of
Lemma 5, it is true in $\ov V$ and $\ov V_\Omega$
that any $\d_\a$ for which $f(\a) = 0$ isn't
$2^{\d_\a} = \d^+_\a$ supercompact.

Let us take this opportunity to observe that the proof of
Theorem 1 uses rather strong hypotheses. Whether a proof of
Theorem 1 is possible from the weaker hypothesis that $\Omega$
is the
least inaccessible limit of supercompact cardinals is
unknown.

\no\S2 The Proof of Theorem 2

We turn now to the proof of Theorem 2. Recall we are assuming
$V \models ``$ZFC + $\k$ is a supercompact limit of
supercompact cardinals + $f : \k \to 2$ is a function''.
As in the remark after Lemma 7 of [A$\infty$a] and
the next to last remark, we also assume without loss of
generality that
$V \models ``$GCH + The supercompact and strongly compact cardinals
coincide except at measurable limit points + Every
supercompact cardinal $\d$ is Laver indestructible under
forcing with $\d$-directed closed partial orderings
not destroying GCH''.
For every $\a < \k$ for which $\a$ is not a measurable
limit of measurable limits of supercompact cardinals, we let
$\d_\a$ and $\la \k^\a_\b : \b < \d_\a \ra$
be as in the proof of Theorem 1.
For every $\a < \k$ for which $\a$ is a measurable
limit of measurable limits of supercompact cardinals,
we let $\d_\a = \a$ but do not define an analogue of
$\la \k^\a_\b : \b < \d_\a \ra$.
$P_\a$ for $\a$ which is not a measurable limit
of measurable limits of supercompact cardinals is then
defined as in the proof of Theorem 1,
and $P_\a$ for $\a$ which is a measurable limit
of measurable limits of supercompact cardinals is
defined as the trivial partial ordering
$\{\emptyset\}$. $P$
is once more defined as the Easton support product
$\underset \a < \k \to{\prod} P_\a$.

\proclaim{Lemma 8}
Let $\a < \k$ be a cardinal which in $V$ is
a regular limit of measurable limits of supercompact cardinals.
Then
$V \models ``\a$ is measurable'' iff
$V^P \models ``\a$ is measurable''.
\endproclaim

\demo{Proof of Lemma 8}
Assume first that $V \models ``\a$ is measurable''.
We show that $V^P \models ``\a$ is measurable''.

Since in $V$, $\a$ is a measurable limit of
measurable limits of supercompact cardinals,
$P_\a$ is trivial. We can thus write
$P_\a = P_{< \a} \times P^\a$.
As $P^\a$ is $\a^+$-directed closed,
$\ov V = V^{P^\a} \models ``\a$ is measurable''.

Let $j : \ov V \to M$ be an elementary embedding
with critical point $\a$ so that
$M \models ``\a$ is not measurable''.
We can then write
$j(P_{< \a}) = P_{< \a} \times Q$,
where the least ordinal $\b_0$ in the domain of $Q$
is so that $\b_0 > \a$. Therefore, if
$H$ is $M$-generic over $Q$ and $G$ is
$\ov V[H]$-generic over $P$,
$j : \ov V \to M$ extends to
$\ov j : \ov V[G] \to M[G \times H]$ in
$\ov V[G \times H]$ via the definition
$\ov j(i_G(\tau)) = i_{G \times H}(j(\tau))$.
We will be done if we can show $H$ is constructible in $\ov V$.

The rest of the argument is similar to the one given in
Lemma 5 of [A$\infty$a]. Specifically, by the fact
GCH holds in $M$ and $M \models ``|Q| = j(\a)$'',
the number of dense open subsets of $Q$ in $M$ is at most
$2^{j(\a)} = {(j(\a))}^+ = j(\a^+)$.
As $\ov V \models \gch$ and $M$ can be assumed to be given
by an ultrapower,
$\ov V \models ``|j(\a^+)| = |{[\a^+]}^\a| = \a^+$''.
Thus, in $\ov V$, we can let $\la D_\g : \g < \a^+ \ra$
enumerate the dense open subsets of $Q$ in $M$.

By the definition of $P_{< \a}$ and the fact
$\b_0 > \a$,
$M \models ``Q$ is $\prec \a^+$-strategically closed''.
As $M^\a \subseteq M$,
$\ov V \models ``Q$ is $\prec \a^+$-strategically closed''
as well. The $\prec \a^+$-strategic closure of $Q$ in both
$\ov V$ and $M$ now allows us to meet all of the dense open
subsets of $Q$ as follows. Work in $\ov V$. Player I picks
$p_\g \in D_\g$ extending
$\sup(\la q_\sigma : \sigma < \g \ra)$
(initially, $q_{-1}$ is the empty condition)
and player II responds by picking
$q_\g \ge p_\g$ (so $q_\g \in D_\g$).
By the $\prec \a^+$-strategic closure of $Q$ in
$\ov V$, player II has a winning strategy for the game, so
$\la q_\g : \g < \a^+ \ra$ can be taken as an
increasing sequence of conditions with $q_\g \in D_\g$
for $\g < \a^+$. Clearly,
$H = \{ p \in Q : \exists \g < \a^+ [q_\g \ge p] \}$
is an $M$-generic object over $Q$ which has been
constructed in $\ov V$.

Assume now that $V^P \models ``\a$ is measurable''.
We show that $V \models ``\a$ is measurable''.
Assume to the contrary that
$V \models ``\a$ is not measurable''.
This implies as earlier in the proof of this lemma
that we can write $P = P_{< \a} \times Q$,
where the least cardinal $\b_0$ in the domain of $Q$
is so that $\b_0 > \a$.
Since $Q$ is therefore $\a^+$-strategically closed
in $V$, GCH and the definition of P imply
$V^P = V^{P_{< \a} \times Q} \models ``\a$ is measurable''
iff
$V^{P_{< \a}} \models ``\a$ is measurable''. Thus, we show
$V^{P_{< \a}} \models ``\a$ is measurable'' implies
$V \models ``\a$ is measurable''.

The argument we use to show
$V^{P_{< \a}} \models ``\a$ is measurable'' implies
$V \models ``\a$ is measurable'' is essentially the
one given in Theorem 2.1.15 of [H] and
Theorem 2.5 of [KiM].
First, note that since $V^{P_{< \a}} \models ``\a$
is Mahlo'',
$V \models ``\a$ is Mahlo''.
Next, let $p \in P_{< \a}$ be so that
$p \force ``\dot \mu$ is a measure over $\a$''.
We show there is some $q \ge p$, $q \in P_{< \a}$
so that for every $X \in
{(\wp(\a))}^V
$,
$q \Vert ``X \in \dot \mu$''.
To do this, we build in $V$ a binary tree
${\cal T}$
of height $\a$, assuming no such $q$ exists.
The root of our tree is $\la p, \a \ra$.
At successor stages $\b + 1$, assuming
$\la r, X \ra$ is on the $\bth$ level of ${\cal T}$,
$r \ge p$, and $X \subseteq \a$, $X \in V$
is so that $r \force ``X \in \dot \mu$'',
we let $X = X_0 \cup X_1$ be such that
$X_0, X_1 \in V$, $X_0 \cap X_1 = \emptyset$,
and for $r_0 \ge r$, $r_1 \ge r$ incompatible,
$r_0 \force ``X_0 \in \dot \mu$'' and
$r_1 \force ``X_1 \in \dot \mu$''.
We can do this by our hypothesis of the non-existence of a
$q \in P_{< \a}$ as mentioned earlier.
We place both $\la r_0, X_0 \ra$ and
$\la r_1, X_1 \ra$ in ${\cal T}$ at height $\b + 1$
as the successors of $\la r, X \ra$.
At limit stages $\l < \a$, for each
branch ${\cal B}$ in ${\cal T}$ of height $\le \l$,
we take the intersection of all second coordinates
of elements along ${\cal B}$. The result is a
partition of $\a$ into $\le 2^\l$ many sets, so since
$\a$ is Mahlo in $V$, $2^\l < \a$, i.e., the
partition is into $< \a$ many sets. Since
$V^{P_{< \a}} \models ``\a$ is measurable'',
there is at least one element $Y$ of this partition
resulting from a branch of height $\l$ and a condition
$s \ge p$ so that $s \force ``Y \in \dot \mu$''.
For all such $Y$, we place a pair of the form
$\la s, Y \ra$ into ${\cal T}$ at level $\l$ as the
successor of each element of the branch generating $Y$.

Work now in $V^{P_{< \a}}$. Since $\a$ is measurable in
$V^{P_{< \a}}$,
$V^{P_{< \a}} \models ``\a$ is weakly compact''.
By construction, ${\cal T}$ is a tree having $\a$
levels so that each level has size $< \a$.
Thus, by the weak compactness of $\a$ in
$V^{P_{< \a}}$, we can let ${\cal B} =
\la \la r_\b, X_\b \ra : \b < \a \ra$
be a branch of height $\a$ through ${\cal T}$.
If we define for $\b < \a$ $Y_\b = X_\b -
X_{\b + 1}$, then since $\la X_\b : \b < \a \ra$
is so that $0 \le \b < \g < \a$ implies
$X_\b \supseteq X_\g$, for $0 \le \b < \g < \a$,
$Y_\b \cap Y_\g = \emptyset$. Since by the construction of
${\cal T}$, at level $\b + 1$, the two second
coordinate portions of the successor of
$\la r_\b, X_\b \ra$ are $X_{\b + 1}$ and $Y_\b$,
for the $s_\b$ so that $\la s_\b, Y_\b \ra$ is at
level $\b + 1$ of ${\cal T}$,
$\la s_\b : \b < \a \ra$ must form in
$V^{P_{< \a}}$ an antichain of size $\a$ in
$P_{< \a}$.

In $V^{P_{< \a}}$, $P_{< \a}$ is a subordering of the
Easton support product
$\underset \b < \a \to{\prod} P_\b$
as calculated in $V^{P_{< \a}}$.
As $V^{P_{< \a}} \models ``\a$ is Mahlo'',
this immediately implies that
$V^{P_{< \a}} \models ``P_{< \a}$ is $\a$-c.c.'',
contradicting that $\la s_\b : \b < \a \ra$ is in
$V^{P_{< \a}}$ an antichain of size $\a$.
Thus, there is some $q \ge p$ so that for every
$X \in {(\wp(\a))}^V$, $q \Vert ``X \in \dot \mu$'',
i.e., $\a$ is measurable in $V$. This contradiction
proves Lemma 8.
\pbf
\hfill $\square$ Lemma 8

By Lemma 8, the measurable limits of $V$-measurable limits
of $V$-supercompact cardinals in $V$ and $V^P$ are
precisely the same. Thus, the proofs of Lemmas 1-6 show
$V^P \models ``$ZFC + If $\a$ is not in $V$ a measurable
limit of measurable limits of supercompact cardinals and
$f(\a) = 0$, then the $\ath$ compact cardinal isn't
supercompact + If $\a$ is not in $V$ a measurable limit
of measurable limits of supercompact cardinals and
$f(\a) = 1$, then the $\ath$ compact cardinal is
supercompact''.

\proclaim{Lemma 9}
$V^P \models ``$Any cardinal $\a \le \k$ which was in
$V$ a supercompact limit of supercompact cardinals is
supercompact''.
\endproclaim

\demo{Proof of Lemma 9}
Since in $V$, $\a$ is a measurable limit of
measurable limits of supercompact cardinals, $P_\a$
is trivial. We can thus write
$P = P_{< \a} \times P^\a$. By the definition
of $P^\a$, the fact all
supercompact cardinals in $V$ are Laver indestructible
under forcing with partial orderings not destroying GCH,
and the fact $P^\a$ is $\a^+$-directed closed,
$V^{P^\a} \models ``\a$ is supercompact''.

Let $\ov V = V^{P^\a}$. The proof of Lemma 9 will be
complete once we have shown
${\ov V}^{P_{< \a}} = V^{P^\a \times P_{< \a}} = V^P
\models ``\a$ is supercompact''. To do this, let
$\l \ge \a$ be arbitrary, and let
$\g = 2^{{[\l]}^{< \a}}$. Let
$j : \ov V \to M$ be an elementary embedding
witnessing the $\g$ supercompactness of $\a$ so that
$M \models ``\a$ is not supercompact''. Note first that
any $\b \in [\a, \g]$ must be so that
$M \models ``\b$ is not supercompact'', for if this were
not the case, then the fact $M^\g \subseteq M$ implies
$M \models ``\a$ is $< \b$ supercompact and $\b$ is
supercompact'' (as $\b$ must be inaccessible in $\ov V$),
so Magidor's theorem of [Ma71] mentioned in
the Introduction tells us
$M \models ``\a$ is supercompact'', a contradiction.
Thus, since $j(P_{< \a}) = P_{< \a} \times Q$, in
$M$, the least cardinal $\b_0$ in the domain of $Q$
must be so that $\b_0 > \g$.

Let $G$ be $V$-generic over $P_{< \a}$ and $H$ be
$\ov V[G]$-generic over $Q$. In
$\ov V[G \times H]$, $j : \ov V \to M$ extends to
$\ov j : \ov V[G] \to M[G \times H]$ via the definition
$\ov j(i_G(\tau)) = i_{G \times H}(j(\tau))$. Since
$M \models ``Q$ is $< \b_0$-strategically closed'' and
$\g < \b_0$, the fact $M^\g \subseteq M$ implies
$\ov V \models ``Q$ is $\g$-strategically closed''
yields that for any cardinal $\sigma \le \g$,
$\ov V[G]$ and $\ov V[G \times H] =
\ov V[H \times G]$ contain the same subsets of $\sigma$.
This means the ultrafilter $\u$ over
${(P_\a(\l))}^{\ov V[G]}$ in $\ov V[G \times H]$
given by $X \in \u$ iff
$\la j(\sigma) : \sigma < \l \ra \in j(X)$
is so that $\u \in \ov V[G]$. This proves Lemma 9.
\pbf
\hfill $\square$ Lemma 9

The proofs of Lemmas 8 and 9 and the remarks
following the proof of Lemma 8 complete the
proof of Theorem 2.
\pbf
\hfill $\square$ Theorem 2

We remark that the proof of Lemma 9 just given requires
no use of GCH. A proof of Lemma 9 using GCH analogous to
the first part of the proof of Lemma 8 can also be given.

\no\S3 Corollaries 3 and 4 and Concluding Remarks

As promised after their statement, we will indicate now
how Corollary 3 using the earlier mentioned weaker
hypotheses and Corollary 4 of Theorem 2 are proven.
Recall that Corollary 3 says from
a measurable limit of measurable limits of supercompact
cardinals, it is consistent that the least measurable
limit of non-supercompact strongly compact cardinals
is the same as the
least measurable limit of supercompact cardinals.
To prove this,
let $V \models ``\k$ is the least measurable
limit of measurable limits of supercompact cardinals''.
Once more, assume without loss of generality that
in addition
$V \models ``\gch$ + The supercompact and strongly compact
cardinals coincide except at measurable limit points +
Every supercompact cardinal $\d$ is Laver indestructible
under forcing with $\d$-directed closed partial orderings
not destroying GCH''. Let $f : \k \to 2$ be given by
$f(\a) = 0$ for even and limit ordinals, and
$f(\a) = 1$ otherwise. Let $P$ be defined as in the proof
of Theorem 2.
Lemmas 1-6 and Lemma 8 then show that $V^P$ is as
desired, with $\k$ by Lemma 8 being the least measurable
limit of both supercompact and non-supercompact
strongly compact cardinals.

To prove Corollary 4, let $V \models ``\k$ is a
supercompact limit of supercompact cardinals'', and once
more, assume the additional hypotheses used in the
proof of Theorem 2.
Let $f : \k \to 2$ be the function which is constantly 0,
and let $P$ be as in the proof of Theorem 2.
If $\k_0$ is in $V^P$ the least supercompact cardinal,
then by the construction of $V^P$,
$V^P \models ``\k_0$ is a limit of strongly compact
cardinals''. This proves Corollary 4.
\pbf
\hfill $\square$ Corollary 3
\pbf
\hfill $\square$ Corollary 4

We note that in the proof of Corollary 4, no use of GCH is
required. The use of GCH in the proof of Theorem 2 is in
the proof of Lemma 8, which in turn is used to show that
if $\k$ is the supercompact cardinal in question, then
the supercompact and strongly compact cardinals below $\k$
satisfy the desired structure properties given by
$f$. If we don't assume GCH
but we assume that
$V \models ``$The supercompact and strongly compact
cardinals coincide except at measurable limit points +
Every supercompact cardinal is Laver indestructible''
and let $f$ be as in the proof
of Corollary 4, since the proof of Lemma 9 requires no use
of GCH, the proof of Corollary 4 just given remains valid.

We take this opportunity to observe that in both Theorems 1
and 2, for any $\a$ so that $f(\a) = 0$, it is possible to
have that the $\ath$ compact cardinal is a bit supercompact
although not fully supercompact. An outline of the
argument for Theorem 1 (we leave it to interested readers
to do the same thing for Theorem 2) is as follows,
assuming we use the notation used in the proof of Theorem 1
and we wish to make the $\ath$ compact cardinal $\d$
when $\d$ isn't supercompact be so
that $\d$ is $\d^+$ supercompact but $\d$ isn't
$\d^{++}$ supercompact:
Let $V \models ``$ZFC + $\Omega$ is the least inaccessible
limit of cardinals $\d$ so that $\d$ is $\d^+$
supercompact and $\d$ is a limit of supercompact cardinals''.
Assume as before without loss of generality that
$V \models ``$GCH + The supercompact and strongly compact
cardinals coincide except at measurable limit points +
Every supercompact cardinal $\d$ is
Laver indestructible under
forcing with $\d$-directed closed partial orderings not
destroying GCH''. Define $P$ as in Theorem 1, except
$\d_\a$ for $\a < \Omega$ is taken as the $\ath$
cardinal $\d$ so that $\d$ is a limit of supercompact
cardinals and $\d$ is $\d^+$ supercompact. The
arguments of Lemmas 1-7, combined with a suitable
generalization of the argument of Lemma 12 of
[AS$\infty$b], will show that $\ov V_\Omega$ is as in
Theorem 1, with the $\ath$ compact cardinal $\d$
being so that if $f(\a) = 0$, then $\d$ isn't
$\d^{++} = 2^{\d^+} = 2^{{[\d^+]}^{< \d}}$
supercompact.

It remains to show that for $\d$ as in the last sentence of
the preceding paragraph, $\d$ is $\d^+$ supercompact in
either $\ov V$ or $\ov V_\Omega$. To see this, we let
$\d =\d_\a$ for some $\a < \Omega$, and we write
$P = P^\a \times P_\a \times P_{< \a}$.
By the amount of strategic closure of $P^\a$, since we
are assuming
$V \models ``\d$ is $\d^+$ supercompact'',
$V^{P^\a} \models ``\d$ is $\d^+$ supercompact''.

An argument analogous to the one found in the first part
of the proof of Lemma 8, with $P_\a$ here taking the place of
the $P_{< \a}$ of Lemma 8, shows
$V^{P^\a \times P_\a} = {(V^0)}^{P_\a} \models
``\d$ is $\d^+$ supercompact''.
If $j : V^0 \to M$ is an elementary embedding witnessing
the $\d^+$ supercompactness of $\d$ so that
$M \models ``\d$ isn't $\d^+$ supercompact'',
$j(P_\a) = P_\a \times Q$, $H$ is $M$-generic over
$Q$, and $G$ is $V[H]$-generic over $P$, then as in
the proof of Lemma 8, $j : V^0 \to M$ extends to
$\ov j : V^0[G] \to M[G \times H]$. We will be done if we can
show $H$ is constructible in $V$, and this is accomplished
via the same sort of argument as in Lemma 8.
Hence,
$V^{P^\a \times P_\a} \models ``\d$ is $\d^+$ supercompact'',
and since $|P_{< \a}| < \d$,
$V^{P^\a \times P_\a \times P_{< \a}} = V^P \models
``\d$ is $\d^+$ supercompact''.

The above paragraph completes our outline.
We leave it to interested readers to fill in any
missing details.

In conclusion, we remark that the proof of Theorem 1
provides a possible plan of attack in obtaining the
relative consistency of the coincidence of the first
$\omega$ measurable and strongly compact cardinals,
or in general, of the relative consistency of the
coincidence of the classes of measurable and
strongly compact cardinals.
If we could show in Lemma 2 that the function $f'$
could be redefined by
$f'(\a) =$ The least measurable cardinal above $\a$
to yield the same sorts of strongly compact
ultrafilters, then if the model $\ov V_\Omega$ of
Theorem 1 were constructed by using
$f : \Omega \to 2$ as the function which is constantly
0 and taking $\la \k^\a_\b : \b < \d_\a \ra$
as the sequence of all measurable cardinals in the interval
$(\underset \b < \a \to{\cup} \d_\b, \d_\a)$,
the model $\ov V_\Omega$ would be so that
$\ov V_\Omega  \models ``$There is a proper class of
measurable cardinals and the classes of measurable
and strongly compact cardinals coincide''.
Of course, the problem of the existence of such
strongly compact ultrafilters is completely open.

\hb\vfill\eject\frenchspacing
\centerline{References}\vskip .5in

\item{[A80]} A. Apter, ``On the Least Strongly Compact
Cardinal'', {\it Israel J. Math 35}, 1980, 225-233.

\item{[A81]} A. Apter, ``Measurability and Degrees of
Strong Compactness'', {\it J. Symbolic Logic 46}, 1981,
180-185.

\item{[A95]} A. Apter, ``On the First $n$ Strongly
Compact Cardinals'', {\it Proceedings Amer. Math. Soc. 123},
1995, 2229-2235.

\item{[A$\infty$a]} A. Apter, ``Laver Indestructibility and
the Class of Compact Cardinals'', to appear in the
{\it J. Symbolic Logic}.

\item{[A$\infty$b]} A. Apter, ``More on the Least Strongly
Compact Cardinal'', to appear in the {\it Mathematical
Logic Quarterly}.

\item{[AG]} A. Apter, M. Gitik, ``The Least Measurable can be
Strongly Compact and Indestructible'', submitted for publication
to the {\it J. Symbolic Logic}.

\item{[AS$\infty$a]} A. Apter, S. Shelah, ``On the Strong Equality
between Supercompactness and Strong Compactness'', to appear in
{\it Transactions Amer. Math. Soc.}

\item{[AS$\infty$b]} A. Apter, S. Shelah, ``Menas' Result is
Best Possible'', to appear in {\it Transactions Amer. Math.
Soc.}

\item{[Bu]} J. Burgess, ``Forcing'', in: J. Barwise, ed.,
{\it Handbook of Mathematical Logic}, North-Holland,
Amsterdam, 1977, 403-452.

\item{[DH]} C. DiPrisco, J. Henle, ``On the Compactness of
$\aleph_1$ and $\aleph_2$'', {\it J. Symbolic Logic 43},
1978, 394-401.

\item{[H]} J. Hamkins, {\it Lifting and Extending Measures;
Fragile Measurability}, Doctoral Dissertation, University of
California, Berkeley, 1994.

\item{[Ka]} A. Kanamori, {\it The Higher Infinite},
Springer-Verlag, New York and Berlin, 1994.

\item{[KaM]} A. Kanamori, M. Magidor, ``The Evolution of
Large Cardinal Axioms in Set Theory'', in:
{\it Lecture Notes in Mathematics 669}, Springer-Verlag,
New York and Berlin, 1978, 99-275.

\item{[KiM]} Y. Kimchi, M. Magidor, ``The Independence
between the Concepts of Compactness and Supercompactness'',
circulated manuscript.

\item{[L]} R. Laver, `` Making the Supercompactness of $\k$
Indestructible under $\k$-Directed Closed Forcing'',
{\it Israel J. Math. 29}, 1978, 385-388.

\item{[LS]} A. L\'evy, R. Solovay, ``Measurable Cardinals and the
Continuum Hypothesis'', {\it Israel J. Math. 5}, 1967, 234-248.

\item{[Ma71]} M. Magidor, ``On the Role of Supercompact and
Extendible Cardinals in Logic'', {\it Israel J. Math. 10},
1971, 147-157.

\item{[Ma76]} M. Magidor, ``How Large is the First
Strongly Compact Cardinal?'', {\it Annals of Math. Logic 10},
1976, 33-57.

\item{[Me]} T. Menas, ``On Strong Compactness and
Supercompactness'', {\it Annals of Math. Logic 7}, 1975,
327-359.

\item{[SRK]} R. Solovay, W. Reinhardt, A. Kanamori,
``Strong Axioms of Infinity and Elementary Embeddings'',
{\it Annals of Math. Logic 13}, 1978, 73-116.
\bye